\documentstyle[11pt,epsfig,rotating,myagsm]{article}

\newcommand{\bZo}{\mbox{\boldmath $Z$}_{\rm obs}}
\newcommand{\bZc}{\mbox{\boldmath $Z$}_{\rm com}}
\newcommand{\obs}{_{\rm obs}}
\newcommand{\com}{_{\rm com}}
\newcommand{\mis}{_{\rm mis}}
\newcommand{\mism}{_{{\rm mis},m}}
\newcommand{\N}{{\cal N}}

\newcommand{\cov}{\mbox{\it Cov}}

\newcommand{\bbf}{\mbox{\boldmath $f$}}

\newcommand{\lm}{\mbox{$m$}}

\newcommand{\bw}{\mbox{\boldmath $w$}}
\newcommand{\bx}{\mbox{\boldmath $x$}}
\newcommand{\by}{\mbox{\boldmath $y$}}

\newcommand{\bI}{\mbox{\boldmath $I$}}

\newcommand{\bW}{\mbox{\boldmath $W$}}
\newcommand{\bR}{\mbox{\boldmath $R$}}

\newcommand{\bSigma}{\mbox{\boldmath $\Sigma$}}

\newcommand{\bdm}{\begin{displaymath}}
\newcommand{\edm}{\end{displaymath}}
\newcommand{\bsy}{\begin{sideways}}
\newcommand{\esy}{\end{sideways}}

\begin{document}
\title{\vspace*{-2cm} Self Consistency: A General Recipe for
  Wavelet Estimation With Irregularly-spaced and/or Incomplete
  Data\thanks{Abbreviated Title: Self Consistent Wavelet Estimation.
  AMS 1991 subject classifications: primary-65T60; secondary-62G08}}
\author{
Thomas C. M. Lee\footnote{Department of Statistics, The Chinese
  University of Hong Kong, Shatin, New Territories, Hong Kong.  Email:
  {\tt tlee@sta.cuhk.edu.hk}}
\\ The Chinese University of Hong Kong/ \\ Colorado State University
\and Xiao-Li Meng\footnote{Department of Statistics, Harvard University,
  Science Center, 1 Oxford Street, Cambridge, MA 02138-2901,
  USA. Email: {\tt meng@stat.harvard.edu}}
\\ Harvard University
}
\date{January 6, 2006}

\maketitle
\vspace*{-0.5cm}
\begin{abstract}
Inspired by the key principle behind the EM algorithm, we propose
a general methodology for conducting wavelet estimation with
irregularly-spaced data by viewing the data as the observed
portion of an augmented regularly-spaced data set. We then invoke
the self-consistency principle to define our wavelet estimators in
the presence of incomplete data.  Major advantages of this
approach include: (i) it can be coupled with almost any wavelet
shrinkage methods, (ii) it can deal with non--Gaussian or
correlated noise, and (iii) it can automatically handle other
kinds of missing or incomplete observations. We also develop a
multiple-imputation algorithm and fast EM-type algorithms for
computing or approximating such estimates. Results from numerical
experiments suggest that our algorithms produce favorite results
when comparing to several common methods, and therefore we hope
these empirical findings would motivate subsequent theoretical
investigations. To illustrate the flexibility of our approach,
examples with Poisson data smoothing and image denoising are also
provided.

\end{abstract}

\section{A Brief Literature Review and Overview}
\label{sec:background} Since the early 1990s, wavelet techniques,
especially for nonparametric regressions and signal denoising,
have attracted enormous attention from researchers across
different fields. Two major reasons for this are that wavelet
estimators enjoy excellent theoretical properties and that they
are capable of adapting simultaneously to spatial and frequency
inhomogeneities (e.g., see
\myncite{Donoho-Johnstone94,Donoho-Johnstone95} and
\myncite{Donoho-et-al95}).  Also,
they are backed up by a fast algorithm (e.g., see
\myncite{Mallat89}).

A restrictive vanilla setting is as follows.  We have $N=2^J$
observations $y_i$ satisfying
\begin{equation}
y_i=f(\frac{i}{N})+e_i, \quad e_i \sim \mbox{ i.i.d.~}\N(0, \sigma^2),
\quad i=0, \ldots, N-1,
\label{eqn:completed}
\end{equation}
and our goal is to estimate the regression function or signal $f$
via a wavelet method. This usually consists of three steps.  The
first step is to compute the empirical wavelet coefficient vector
$\bw$ by applying a discrete wavelet transform (DWT) to $\by=(y_0,
\ldots, y_{N-1})^T$; that is, if $\bW$ denotes the DWT matrix, then
$\bw$ is given by $\bw=\bW\by$.  The second step is to apply a
shrinkage operation (e.g., thresholding) to $\bw$ to obtain an
estimated wavelet coefficient vector $\hat{\bw}$. Finally,
$\hat{\bbf}=\bW^T\hat{\bw}$ is computed as the wavelet estimate or
reconstruction of $\bbf=(f(0), \ldots, f(\frac{N-1}{N}))^T$.
Details can be found in numerous monographs, such as
\citeasnoun{Daubechies92},
\citeasnoun{Ogden96}, \citeasnoun{Mallet99},
and \citeasnoun{Vidakovic99}.

The most popular shrinkage operation in the statistical literature
appears to be thresholding.  Earliest examples include the
``universal'' thresholding method of
\citeasnoun{Donoho-Johnstone94}, the SURE thresholding method of
\citeasnoun{Donoho-Johnstone95}, and the method of
\citeasnoun{Saito94} that uses the minimum description length
(MDL) principle of \citeasnoun{Rissanen89}.  Since then many
different thresholding methods have also been proposed, including
the cross--validation method of \citeasnoun{Nason96}, the refined
MDL based methods of \citeasnoun{Moulin96} and of
\citeasnoun{Lee02:wavetree}, the cross--validatory {\sc AIC}
method of \citeasnoun{Hurvich-Tsai98}, and the method of
\citeasnoun{Crouse-et-al98} that uses hidden Markov models.
Moreover, various Bayesian and empirical Bayes methods have also been
proposed; e.g., see \citeasnoun{Chipman-et-al97},
\citeasnoun{Abramovich-et-al98}, \citeasnoun{Clyde-et-al98},
\citeasnoun{Vidakovic98}, Clyde \& George (1999,
2000)\nocite{Clyde-George99,Clyde-George00}, and
\citeasnoun{Johnstone-Silverman05}. 
In addition, some treatment has also been given to the issue of
correlated noise, see, for example, \citeasnoun{Wang96},
\citeasnoun{Johnstone-Silverman97}, \citeasnoun{Jansen-Bultheel99}
and \citeasnoun{Lee02:wavetree}.
Lastly, robust wavelet smoothing has been considered for examples by
\citeasnoun{Sardy-et-al01} and \citeasnoun{Oh-et-al05},
while methods for reducing boundary artifacts are studied by
\citeasnoun{Lee-Oh04} and \citeasnoun{Oh-Lee05}.

However, the applicability of many existing wavelet regression
methods is limited by the assumptions that the data are not only
complete but also equispaced and that the number of design points
is an integer power of 2, as in (\ref{eqn:completed}).  Different
methods have been proposed to relax the equal-spaced assumption.
\citeasnoun{Hall-Turlach97} proposed the uses of two interpolation
rules for mapping the observed data into a regular grid.  They
provided a detailed theoretical analysis of their proposals and an
asymptotic choice of the thresholding value.
\citeasnoun{Kovac-Silverman00} also investigated the use of
interpolation.  They developed a fast algorithm for computing the
noise covariance structure after the original observed data are
mapped into an interpolation grid. With the knowledge of this new
noise covariance structure, tailored thresholding values can be
derived. \citeasnoun{Nason02} demonstrated that this fast
covariance updating algorithm is extremely useful for speeding up
cross--validation type calculations.  Interpolation based methods
were also studied in \citeasnoun{Sardy-et-al99}.
\citeasnoun{Cai-Brown98} took a different approach by invoking
distributional assumptions for the design points. More recently,
\citeasnoun{Antoniadis-Fan01} consider the use of penalized least
squares for handling non--equally spaced data.  In their procedure
the best curve estimate is defined as the minimizer of a
regularized least squares criterion.  Starting with a regularized
wavelet interpolation of the non--equally spaced observed data, a
one--step iterative algorithm is applied to approximate the
minimizer. Finally, the method of lifting can be applied to
construct wavelets for irregularly spaced data (e.g., see
\myncite{Delouille-et-al04} and \myncite{Nunes-et-al06}). However,
thresholding methods for such lifting wavelet bases seem to be less
developed.

Whereas most of these non-equispaced methods are effective for
various applications, they impose some additional assumptions that
are not used with regular designs; e.g., the implicit smoothness
assumption in interpolation methods.  In this article we
demonstrate that it is possible to deal with the irregular design
problem without making such assumptions. Our key idea is to view
an irregular design as a regular design with missing observations.
This view allows us to utilize those well--known and widely tested
methods in the extensive literature on estimation and computation
for missing data, in particular EM--type algorithms (e.g.,
\myncite{demp:lair:rubi:77,meng:rubi:93} and \myncite{meng:vand:97})
and multiple imputation (e.g.,
\myncite{rubi:87} and \myncite{meng:94}).  We invoke a
self--consistency criterion, in the same spirit as in
\citeasnoun{Efron67} and as the self--consistency principle
underlying the EM algorithm, to {\it define} the wavelet
regression estimator with incomplete data.

Our approach is that, given only a complete--data procedure and
the corresponding model on the missing--data mechanism, we first
seek the most efficient wavelet estimator for the values of the
regression function at the observed designed points; i.e., design
points at which the response $y$ is observed.  We then incorporate
into such ``optimal'' estimation procedure any additional
assumptions that are not used with the complete procedure,
provided that such assumptions can further improve the efficiency
for any particular applications. This {\it separation} between the
inherited information in the observed data and the information
built into a procedure due to external assumptions helps to obtain
more efficient wavelet regression estimates with irregular
designs.  This is not surprising as sensible self--consistent
procedures often lead to the most efficient estimators,  both in
the parametric estimation (e.g., maximum likelihood estimation via
the EM algorithm) and the nonparametric estimation (e.g., the
Kaplan-Meier estimator) contexts.  Indeed, as
\citeasnoun{Tarpey-Flury96} put it, self--consistency is a
fundamental concept in statistics, and is a general statistical
principle for retaining as much as possible the information in the
data.  In addition, as demonstrated below, another advantage of
this approach is that it can be straightforwardly extended to more
complicated settings, such as image denoising and non-Gaussian
errors, and it handles non-equispaced data simultaneously with
other kind of incomplete data, such as in photo inpainting
applications (see Section~\ref{sec:sim}).

The rest of this article is organized as follows.
Section~\ref{sec:SC} reviews the self--consistency principle.  It also
presents a self--consistency criterion for regression function
estimation, from which our non-equispaced wavelet
estimator is implicitly defined.  In Section~\ref{sec:alg} three
algorithms are constructed to compute or approximate this wavelet
estimator.  Further extensions including two-dimensional
implementation are discussed in Section~\ref{sec:extensions}, and
simulation results are reported in Section~\ref{sec:sim}.  Future
work, especially regarding theoretical development, is discussed in
Section~\ref{sec:conclude}.

\section{Self Consistency: How Does It Work?}
\label{sec:SC}
\subsection{Self-consistency: An Intuitive Principle}

To illustrate the self-consistency principle, imagine the following
scenario.  Mr.~Littlestat was told by his boss to prepare a
presentation on the growth in sales since the company's inception in
1993, and to make a prediction for the next couple of years.
Eye-balling the 13 years of the data he was given, he felt that he
could draw a reasonably-looking line, but he also knew that it would
not please his boss.  He vaguely remembered something called
``least-squares line'' from his college days, but that was all he
could remember.  So he asked his teenage son, who was a member of a
high-school math club.  Indeed, his son not only knew about the
method, but actually had just programmed it for a homework problem.
However, as his son was dashing out for a movie date, he briefly
showed the program on his laptop to his father and said ``Dad, I'm
really running late.  Just type in your data as two columns here,
click that little thing, and you will find a drawing at the printer.''

Not thrilled by his son's rushing but nevertheless pleased to have
the program, Mr.~Littlestat sat down and started to type in his
data.  Then he was really unhappy, as he found out that his son's
program was hardwired specifically for solving the homework that had
16 data points.  It simply would not run for Mr.~Littlestat data set
with 13 data points.  There were three more rows needed to be filled
in.

``What should I do now?''  Mr.~Littlestat asked himself.  He knew
nothing about the method, nor how to modify his son's program.  Nor
did he have the patience to wait for his son's return, as he really
needed to finish the preparation for his presentation tomorrow.

As the desperation set in, Mr.~Littlestat thought ``Well, what if
I just make up some numbers for the next three years, and see what
happens?''  So he did, and clicked.  Instantly, a printout came out
from the printer with a drawing of a line and the 16 points displayed,
including the three fake sales figures.

Excited, Mr.~Littlestat examined the plot, and saw the line was
visually a bad fit to the 13 years of the real data. ``Well,
that's expected, since I put in some fake sales figures'', he murmured
to himself.  But then it hit him: ``since I made up these three points
anyway, why don't I just make them to sit on this line and run the
program again? That probably would be better ...''

So he reentered the three fake sales figures by reading off the sales
figures from the line, clicked again. The new printout showed a line
fitting better for the past 13 years, but the three fake sales were
still off from the new line, though they were a bit closer. ``Hmmm,
this is interesting''.  He was getting intrigued by his clever
invention, ``why not try this again?''  So he reentered the three
future sales by reading off the line again.

Mr.~Littlestat's excitement increased with the number of printouts, as
the line fluctuated less and less, and the three future sales got
closer and closer to the fitted line.  Then everything stopped, and
the three future figures sit on the line exactly, as far as
Mr.~Littlestat could tell.  Visually inspecting all the plots he had,
he exclaimed, ``I guess this is it!'' while holding the last plot.

Although Mr.~Littlestat was far from sure whether the line in his
hand had anything to do with the least-squares line he was after,
he convinced himself that he had done his best given what he was
given, because there was really nothing else he could do --- the
line simply stopped moving, and that had to be its best in some sense
for it could not be improved further.

Mr.~Littlestat was indeed right.  His intuitive pursue actually had
led to the correct answer.  The limit of his iterative procedure is
indeed the least-squares fit, a consequence of the
``self-consistency'' property, as we shall explain below.  We provided
this story to illustrate and emphasize that self-consistency is
nothing more than good common sense.  It does not always work, just as
not all ``common senses'' would lead to good answers, but it is
typically suggestive, and often leads to optimal solutions that can be
justified mathematically.

In the least-squares example above, Mr.~Littlestat's method worked
because the least-squares estimator is self-consistent in the
following sense.  Suppose we have a regression setting for which $x$
is univariate:
\[
y_i=\beta x_i + \epsilon_i, \quad i=1, \ldots, n,  \quad
\epsilon_i \sim \mbox{ i.i.d.~}F(0, \sigma^2),
\]
where $F(0,\sigma^2)$ denotes a distribution with mean zero and variance
$\sigma^2$. The least-squares estimator of $\beta$ then is given by
\begin{equation}\label{eqn:leas}
\hat\beta_n\equiv \hat\beta_{n}(y_1, \ldots, y_{n}) =
\frac{\sum_{i=1}^{n}y_ix_i} {\sum_{i=1}^{n} x_i^2}.
\end{equation}
Then, for any $m< n$, as long as $\sum_{i=m+1}^{n} x_i^2 >0$,
\begin{equation}\label{eqn:selfl}
E\left(\hat\beta_{n}\Big\vert y_1, \ldots, y_{m};
\beta=\hat\beta_{m} \right) = \hat\beta_{m}.
\end{equation}
That is, the least-squares estimator has a Martingale-like
property, and reaches a perfect equilibrium in its projective
properties. Therefore, we can obtain $\hat\beta_m$ from a
procedure for computing $\hat\beta_n$ with $n>m$ by
solving~(\ref{eqn:selfl}) for the ``fixed point'' $\hat\beta_m$.
And this can be solved iteratively {\it without knowing} the form
of $\hat\beta_{n}$ as long as we can compute the average of the
{\it values} of $\hat\beta_{n}$ over the conditional distribution
as required by the left-hand side of~(\ref{eqn:selfl}).  In this
specific case, Mr.~Littlestat's ``line fitting'' method is correct
because of the linearity of $\hat\beta_n$ in the $y_i$'s.
Specifically, starting with some initial guess, $\beta_{13}^{(0)}$
say, at the $t^{th}$ iteration, we can impute the three
``missing'' $y_i$'s  by their conditional expectations
$y_i^{(t)}=\beta_{13}^{(t)}x_i$ ($i=14, 15, 16$), and then compute
the next iterative estimate
\begin{equation}\label{eqn:iter}
\beta_{13}^{(t+1)} = \hat\beta_{16}(y_1, \ldots, y_{13},
y_{14}^{(t)}, y_{15}^{(t)}, y_{16}^{(t)}).
\end{equation}
Evidently, because of (\ref{eqn:leas}), the limit of
(\ref{eqn:iter}), denoted by $\hat\beta_{13}$, must satisfies
$$ \hat\beta_{13} =  \frac{\sum_{i=1}^{13}
y_ix_i + \hat\beta_{13} \sum_{i=14}^{16} x_i^2}{\sum_{i=1}^{16}
x_i^2},$$ which means that
$\hat\beta_{13}=\sum_{i=1}^{13}y_ix_i/(\sum_{i=1}^{13}x_i^2)$,
the correct least-squares estimate with 13 data points.

\subsection{A Self--consistent Regression Estimator}
Inspired by the least-squares estimator above, we propose to also
invoke the self--consistency principle for regression function
estimation in general when facing incomplete data. Note that
although in this paper we focus on wavelet regression models, our
proposal extends to more general non-parametric or semi-parametric
regression setting. Specifically, let $\bZo$ and $\bZc$ be,
respectively, the observed and (imaginary) complete data. Suppose
for the moment that, given $\bZc$, we have a method for computing
the ``best'' estimate $\hat f\com$ for our regression function
$f$. We propose that $\hat f\obs$, our estimate of $f$ given
$\bZo$, under squared loss, to be the solution of the following
self--consistent equation:
\begin{equation}
E\left\{ \hat f\com(\cdot)\big\vert \bZo,f=\hat f\obs,
\theta=\hat\theta\obs \right\} = \hat f\obs (\cdot),
\label{eqn:selfcon}
\end{equation}
where $\theta$ collects all nuisance parameters (such as the
variance parameter); for notation simplicity, for the rest of this
paper we will suppress, but not ignore, this conditioning on
$\theta$.  Since (\ref{eqn:selfcon}), as demonstrated in
Section~\ref{sec:alg}, can be solved numerically via iterations,
it provides a way to obtain our ``best'' incomplete-data estimator
$\hat f\obs$ by simply using the corresponding complete-data
procedure that computes $\hat f\com$, much like the EM algorithm
obtains the incomplete-data MLE via procedures for complete-data
MLE. In doing so, no additional assumptions will be required,
other then the necessary specification of the conditional
distribution of the missing part of $\bZc$ given the observed
$\bZo$.  We note that such specifications are necessary as
otherwise what is missing can never be recovered.

The self--consistent equation (\ref{eqn:selfcon}) was also
motivated by its success in estimating cumulative distribution
function (CDF) with censored or truncated data, where
$f\com(\cdot)$ will be the empirical CDF, a topic that has been
studied extensively in the literature (e.g., see
\myncite{Efron67}).  Indeed, it is well known that the
Kaplan--Meier estimator is the solution to (\ref{eqn:selfcon}) for
the right censored data under non-informative censoring, and it is
also the generalized maximal likelihood estimator (e.g.,
\myncite{joha:78}). In addition, similar ideas have also been
successfully adopted to construct various nonlinear data
summaries; e.g., see \citeasnoun{Tarpey-Flury96} and references given
therein.  Of course, the most spectacular success of the
self--consistency principle is the EM algorithm
\cite{demp:lair:rubi:77} and its various generalizations (e.g.,
\myncite{meng:rubi:93} and \myncite{meng:vand:97})
where the self-consistency principle, when applied to the
complete-data score function, typically leads to incomplete--data
maximum likelihood estimator. Indeed, our whole investigation as
reported in this paper was guided closely by our knowledge and
insights of the development of EM algorithm and its extensions,
which turned out to be extremely fruitful for our current
purposes.

\subsection{Heuristics}

In view of these great successes, it is natural to expect that our
estimator $\hat f\obs$, the solution to (\ref{eqn:selfcon}), would
possess excellent statistical properties.  The following heuristics
provides some insight and indication on why (\ref{eqn:selfcon}) can
lead to efficient estimator under squared loss.  Suppose $\hat f\com$
is the complete-data optimal estimator that minimizes $\|f-\hat
f\com\|^2$. Then, for any $\hat f\obs$,
\begin{eqnarray*}
\|f-\hat f\obs\|^2
& = & \|f-\hat f\com\|^2 + \|\hat f\com - \hat f\obs\|^2 \\
& = & \|f-\hat f\com\|^2 + \|\hat f\com - E(\hat f\com|\bZo, f)\|^2 +
\|E(\hat f\com|\bZo, f) - \hat f\obs \|^2.
\end{eqnarray*}
Thus, minimizing $\|f-\hat f\obs\|^2$ over $\hat f\obs$ is
equivalent to minimizing $\|E(\hat f\com|\bZo, f) - \hat
f\obs\|^2$, because the first two terms in the right hand side of
the last equality are constants with respect to the minimization.
It is thus natural to suspect that the solution of
(\ref{eqn:selfcon}) is the minimizer of $\|f-\hat f\obs\|^2$, at
least asymptotically.

There is also a Bayesian heuristics for (\ref{eqn:selfcon}).  From
a Bayesian view point, if $\hat f\com$ is the Bayesian estimator
for $f$ given the complete data, then the Bayesian estimator given
the observed data, under the squared loss, is $\hat
f\obs(\cdot)=E\{\hat f\com(\cdot)|\bZo\}$.  Although $\hat f\obs$
depends on the choice of prior, if it is an efficient estimate of
the true $f$, it is reasonable to expect that asymptotically the
ratio between the posterior expectation and the conditional
sampling expectation evaluated at $f=\hat f\obs$, both of $\hat
f\com$, converges (almost surely) to~1. That is,
\begin{equation}
\frac{E\left\{\hat f\com(\cdot)\big\vert \bZo, f=\hat f\obs
    \right\}}{\hat f\obs(\cdot)}
=\frac{E\left\{\hat f\com(\cdot)\big\vert \bZo, f=\hat f\obs \right\}}
{E\left\{\hat f\com(\cdot)\big\vert \bZo\right\}}
\rightarrow 1,
\label{eqn:selrati}
\end{equation}
which implies that (\ref{eqn:selfcon}) should hold at least
asymptotically.

We note that although the heuristics arguments above are
appealing, we currently do not have rigorous theoretical results
to establish the optimality of the self-consistent estimator.
Indeed, even to prove the existence of the solution to the
self-consistent equation (\ref{eqn:selfcon}) is a theoretical task
that we have not been able to carry through. Our work, as reported
in this paper, therefore has been more of an ``engineering
nature'', focusing on constructing algorithms to solve
(\ref{eqn:selfcon}) and to demonstrate empirically the good
performance, conceptual and implemental simplicity, as well as
the flexibility of the self-consistent approach. We hope that
these empirical demonstrations show the great promise of this
self-consistent approach, and thereby stimulate investigation of
the theoretical properties, including optimality, of the
self-consistent estimator.

\section{Three Algorithms}
\label{sec:alg} The self-consistent estimator $\hat f\obs$ would
not be of much practical value if (\ref{eqn:selfcon}) could not be
solved with reasonable computational effort.  To solve
(\ref{eqn:selfcon}), two steps are involved.  The first is to
carry out the conditional expectation on the left-hand side, and
the second is to solve the equation, in analogous to the E-step
and M-step of the EM algorithm, respectively.  However, unlike
many common EM applications where the E-step is in closed form, in
the wavelet applications, the exact E-step is typically
analytically infeasible.  This is because, due to the shrinkage
operation, $\hat f\com$ is a highly complicated non-linear
function of the missing $y's$.  There are two general approaches
for dealing with such a problem.  The first is to use a Monte
Carlo E-step, as in \citeasnoun{wei:tann:90} and
\citeasnoun{meng:schi:96}, and the second is to trade the
exactness for simplicity by making certain approximations to the
conditional expectation.  Below we propose three algorithms for
computing or approximating $\hat f\obs$, one of which is based on
the first Monte Carlo approach, and the other two follow the
second approximation approach.

\subsection{A Multiple Imputation Self-Consistent (MISC) Algorithm}
\label{sec:MISC} First we fix the notation.  Define $\bI\obs$ as
the ``observed data index set'': $i\in \bI\obs$ if the $i$th data
point $(x_i, y_i)$ in~(\ref{eqn:completed}) is observed.  Let
$\by=(y_0, \ldots, y_{N-1})^T$ denote the complete responses, and
let $\by\mis$ and $\by\obs$ denote respectively the missing and
observed portions of $\by$.  That is, $\by\mis=\{y_i: i\not\in
\bI\obs\}$ and $\by\obs=\{y_i:i\in\bI\obs\}$.  Define $\bx\obs =
\{x_i : i\in \bI\obs\}$, and hence the observed data is
$\bZo=\{\bx\obs,\by\obs\}$. Also denote the covariance matrix of
$\by$ as $\cov(\by)=\bSigma$.  Thus we allow the error terms
$e_i$'s to be correlated, as long as $\bSigma$ can be identified
and efficiently estimated from $\by\obs$.

Our first algorithm, termed the {\em multiple imputation
  self-consistent (MISC) algorithm}, assumes that a
complete-data wavelet regression
procedure has been chosen (e.g., the SURE method of
\myncite{Donoho-Johnstone95}).  Starting with initial estimates $\hat
f^{(0)}$ and $\hat{\bSigma}^{(0)}$, the algorithm iterates the
following three steps for $t=1, \ldots$:
\begin{description}
\item[Step~1] {{\sl Multiple Imputation}:} For $\lm=1, \ldots, M$,
simulate $\by\mism$ independently from
\[P(\by\mism\big\vert \by\obs; f=\hat f^{(t-1)},
\bSigma=\hat{\bSigma}^{(t-1)}).\]

\item[Step~2] {\sl Wavelet Shrinkage:} For $\lm=1, \ldots, M$,
apply the chosen complete-data wavelet shrinkage procedure to the
{\em completed data} $\by_{\lm}=\{\by\obs, \by\mism\}$ and obtain
$\hat f_{\lm}(x_i)$, $i=0, \ldots, N-1$.

\item[Step~3]{\sl Combining Estimates:}
Compute the $t$-th iterative estimate of $f$ as
\begin{equation}
\hat f^{(t)}(x_i) = \frac{1}{M} \sum_{\lm=1}^M \hat f_{\lm}(x_i),
\quad i=0, \ldots, N-1.
\label{eqn:combine}
\end{equation}
Also, use the residuals $\{ y_i - \hat f^{(t)}(x_i):
i\in\bI\obs\}$ to obtain an efficient estimate
$\hat{\bSigma}^{(t)}$, such as MLE, of $\bSigma$.
\end{description}
In Step~1 above, the larger the $M$, the better results one would
expect, but at the expense of increased computational time.  Our
numerical experience indicates that, as long as $M$ is larger
than a minimum cutoff, the additional improvement on $\hat f$
computed with a larger $M$ is not largely significant. In our
numerical experiments, we typically used $M=100$, though $M=10$ is
sometimes acceptable as well.

It is evident that the above is a generic algorithm, in the sense
that it is not restricted by the specific form of the
complete-data wavelet regression procedure, nor by the choice of
the error distribution.  On one hand, this is a great advantage as
it is extremely flexible and the additional programming, relative
to that for the complete-data procedure, is minimal as long as it
is easy to draw from the conditional distribution in the first
step (which typically is the case for Gaussian errors, independent
or not). It also provides a benchmark and basis for developing
more specialized and sophisticated algorithms.  On the other hand,
it is a ``brute force'' algorithm, and is thus quite inefficient
as a numerical algorithm.

\subsection{A Simple (Sim) Approximated Algorithm}
\label{sec:Sim}
To construct a faster algorithm for computing
$\hat f\obs$, we started with replacing the costly multiple
imputation step in the MISC algorithm by a very simple analytical
approximation. We label the resulting algorithm as the {\em simple
(Sim) algorithm}, and it is designed for a specific type of
shrinkage methods, namely, for hard thresholding methods for which
the thresholding value is a known function $g(\hat\sigma)$ of
$\hat\sigma$, where $\hat\sigma$ is an estimate of $\sigma$.  A
classical example for $g(\hat\sigma)$ is the universal
thresholding scheme of \citeasnoun{Donoho-Johnstone94}, for which
$g(\hat\sigma) = \hat\sigma \sqrt{2\log N}$.

Starting with $\hat{f}^{(0)}$ and $\hat{\sigma}^{(0)}$, the Sim
algorithm iterates, at the $t$-th iteration, the following steps:
\begin{description}
\item[Step~1] For each $i$ such that $i\not\in\bI\obs$, impute the
  corresponding missing $y_i$ by $y_i^{(t)}=\hat f^{(t-1)}(x_i)$,
  which creates a completed-data set: $\by^{(t)}=\{y_i: i\in
  \bI\obs\}\cup\{y_i^{(t)}: i\not\in \bI\obs\}.$

\item[Step~2] Apply a DWT to $\by^{(t)}$ to obtain the empirical
  wavelet coefficients $\bw^{(t)}=\bW\by^{(t)}$.

\item[Step~3] Obtain a robust estimate $\tilde\sigma^{(t)}$ of
  $\sigma$ from $\bw^{(t)}$, for example, the median absolute
  deviation method used by \citeasnoun{Donoho-Johnstone94}. We call
  $\tilde\sigma^{(t)}$ the {\em unadjusted} estimate for $\sigma$.

\item[Step~4] Use the following {\em variance inflation formula} to
  obtain an {\em adjusted} estimate $\hat\sigma^{(t)}$ for $\sigma$:
\begin{equation}
\hat\sigma^{(t)} = \sqrt{\{\tilde\sigma^{(t)}\}^2 +
C_m \{\hat\sigma^{(t-1)}\}^2},
\label{eqn:inflate}
\end{equation}
where $C_m = 1-\frac{n}{N}$ is the fraction of missing observations.

\item[Step~5] Compute $\hat{\bw}^{(t)}$ by thresholding $\bw^{(t)}$
  with the thresholding value $g(\hat\sigma^{(t)})$.

\item[Step~6] Apply the inverse DWT to $\hat{\bw}^{(t)}$ and obtain
  the $t$-th iterative estimate
  $\hat{\bbf}^{(t)}=\bW^T\hat{\bw}^{(t)}$.
\end{description}
In all our numerical experiments, convergence was declared if
$|\hat\sigma^{(t)}-\hat\sigma^{(t-1)}|/\hat\sigma^{(t)}<\epsilon$.
Upon convergence, estimates of $f$, as well as $\sigma$, will be
obtained. It is obvious that computationally this Sim algorithm is
much less intensive than MISC, because it only requires one
complete-data wavelet shrinkage computation within each iteration,
in contrast to the $M$ sets of computation required by MISC.

A key component of the Sim algorithm is the variance inflation
formula (\ref{eqn:inflate}), which takes into account the effect
of those imputed $y_i^{(t)}$'s on the estimation of $\sigma^2$.
The formula was borrowed from the EM algorithm for estimating
$\sigma^2$ with normal regression under independent errors,
which requires replacing each missing $y_i^2$ by its conditional
expectation
\[
E\left[y^2_i\bigg|\by\obs; f=\hat f^{(t-1)},
\sigma^2=\{\sigma^{(t-1)}\}^2 \right]=\{y_i^{(t)}\}^2 +
\{\sigma^{(t-1)}\}^2.
\]
Replacing the complete-data sufficient statistics $\sum_i y_i$ and
$\sum_i y_i^2$ in the complete-data MLE for $\sigma^2$ then lead
to (\ref{eqn:inflate}). Although we recognize that this adjustment
may not be consistent with the method used for obtaining the
unadjusted $\tilde\sigma^{(t)}$, which often is not MLE, we adopt
(\ref{eqn:inflate}) for its conceptual and implemental
simplicity. As demonstrated in Section~\ref{sec:sim}, it works
quite well in the sense that not carrying out this variance
inflation adjustment would lead to noticeably poorer wavelet
estimates. However, this variance inflation adjustment does not
account for all the uncertainty in the thresholding due to
missing data, and hence it does not work well when the percentage
of missing data is large. A better approximation therefore is
needed.

\subsection{A Refined (Ref) Fast Algorithm}
\label{sec:Ref}

To simplify presentation, we will use single-indexing $w_l$
instead of the usual double-indexing $w_{jk}$ notation to denote a
wavelet coefficient. At the $t$-th step, we need to calculate
\begin{displaymath}
\hat f^{(t)} = E\left\{\hat f\com \Big| \by\obs, f=\hat
  f^{(t-1)}\right\},
\end{displaymath}
which amounts to calculating all
\begin{equation}\label{eqn:estep}
\tilde w_{l}^{(t)} \equiv E\left\{1_{|w_{l}|\ge
g(\tilde\sigma)}w_{l}
  \Big| \by\obs, f=\hat f^{(t-1)}\right\},
\end{equation}
where $w_{l}$'s and $\tilde\sigma$ are respectively the
complete-data empirical wavelet coefficients and variance
estimate. The simple algorithm in Section~\ref{sec:Sim} took a
very crude approximation of this conditional expectation, by
thresholding the conditional expectation of $w_{l}$ with
$g(\hat\sigma)$ using the adjusted $\hat\sigma$. We can obtain
better approximations if we are willing to give up some generality
(but see Section~\ref{sec:non}).

For example, under the i.i.d.~normal error setting of
(\ref{eqn:completed}), we can obtain a much refined analytic
approximation if we ignore the conditional variability in $\tilde
\sigma$ when calculating the conditional expectation for
(\ref{eqn:estep}). To proceed, we need to define the quantity
$\eta_l$: if $w_l^{(t)}$ is the empirical wavelet coefficient
obtained from Step~2 of the Sim algorithm, then the conditional
distribution of $w_l$ given $\{\by\obs, \sigma^2\}$ is
$\N(w_l^{(t)}, \eta_l^2\sigma^2)$. With this setup and denote
$g(\sigma)$ by a constant $c$ to signify the fact that we assume
it is known, we shown in Appendix~\ref{sec:closed} that
(\ref{eqn:estep}) has the following simple analytic expression
\begin{equation}\label{eqn:ce1}
\tilde w_{l}^{(t)} = \alpha(w_{l}^{(t)}, \eta_l) + \beta(w_{l}^{(t)},
\eta_l) \times w_{l}^{(t)},
\end{equation}
where the $\alpha$ and $\beta$ functions are given by
\begin{equation}\label{eqn:ce2}
\alpha(w, \eta) = {\eta\sigma\over \sqrt{2\pi}}
\left\{e^{-{1\over 2}\left({c-w\over \eta\sigma}\right)^2}
- e^{-{1\over 2} \left({c+w\over \eta\sigma}\right)^2}\right\}
\quad {\rm and} \quad
\beta(w, \eta) = 2 - \Phi\left({c-w\over \eta\sigma}\right) -
\Phi\left({c+w\over \eta\sigma}\right),
\end{equation}
with $\Phi$ being the CDF function of $\N(0,1)$.

The resulting algorithm is identical to Sim except that we replace
its Step~5 by (\ref{eqn:ce1}) and (\ref{eqn:ce2}), where we use
$\sigma=\hat\sigma^{(t)}$ and $c=g(\hat\sigma^{(t)})$. Thus,
computationally, this new refined (Ref) algorithm is almost as
efficient as Sim, and it is also straightforward to program as
only standard functions are involved in (\ref{eqn:ce1}) and
(\ref{eqn:ce2}). However, because it provides a much more refined
E-step, the statistical efficiency of the resulting estimator is
expected to be much closer to that of the MISC estimator with
$M=\infty$.  Numerical results from Section~\ref{sec:sim} strongly
support this expectation.

The analytic formula (\ref{eqn:ce1})-(\ref{eqn:ce2}) deserves
several important remarks. First, under the assumption of i.i.d.~noise
error, the value of $\eta_l$ can be easily computed at the
outset of iteration as the $l$th diagonal element of $\bI-\bW \bR
\bW^\top$, where $\bR$ is an $N\times N$ matrix whose off-diagonal
elements are all zero, and whose $l$th diagonal elements is one if
$y_l$ is observed and zero otherwise; that is, the diagonal of
$\bR$ forms a ``response indicator" vector for $\by$. Note that
the assumption of i.i.d.~error is non-essential as long as the
error covariance $\bSigma$ is known, which we do assume for the
current approximation. Intriguingly, approximating all $\eta_l$'s
by their average, which is exactly the fraction of missing data
$C_m$ as used in the variance inflation formula
(\ref{eqn:inflate}), works surprisingly well -- see
Section~\ref{sec:2d} for more discussion.

\begin{figure}[ht]
\begin{center}
\vspace*{-0.7cm}
\epsfig{file=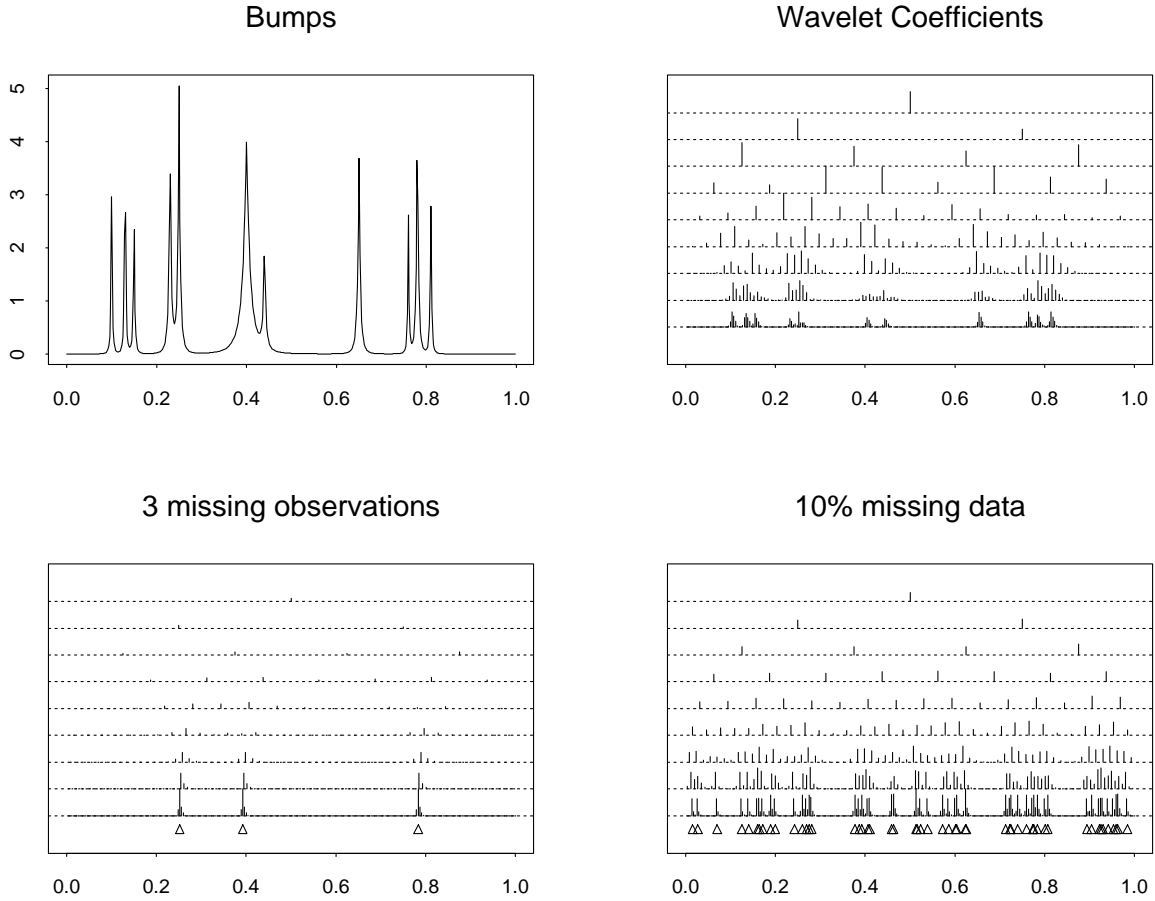,width=6.5in}
\vspace*{-2cm}
\end{center}
\caption{\small Frequency--location plots of $\sqrt{\eta_l}$.  Top
left: test function {\em Bumps}; top right: wavelet coefficients
of {\em Bumps}; bottom left: plots of $\sqrt{\eta_l}$ with three
missing observations; bottom right: plots of $\sqrt{\eta_l}$ with
10\% missing observations. In the $\eta_l$ plots those triangles
at the bottom indicate locations of missing observations.  Note
that $\sqrt{\eta_l}$ is plotted instead of $\eta_l$ to enhance
visibility.} \label{fig:eta}
\end{figure}

Second, as a by-product, the quantity $\eta_l$ can be seen as a
measure of the percentage of missing information in $w_l$ due to
missing data.  It is because $0\le \eta_l \le 1$, and that $\eta_l$ is
one when there is no information in the observed data about $w_l$ and
zero when $w_l$ is fully observed.  Since the missing data here are
the result of irregular design points, we can also view $\eta_l$ as a
measure of {\em irregularity} in the data for the wavelet coefficient
$w_l$. Figure~\ref{fig:eta} provides some illustrative
plots of $\eta_l$, and demonstrates well that the effects of the
missing data, as expected, are localized. Also, the plots indicate
that these missing data seem to have a stronger impact on high
frequency wavelet coefficients. We believe that this interesting
by-product, that is, the ``irregularity plot" (i.e., by displaying
the $\eta_l$'s on a frequency-location plot, as we do with the
empirical wavelet coefficients), is worthy of further exploring,
for example, for the purpose of diagnosing and determining the
suitability of a particular wavelet regression model for a
particular irregular design.

Third, an intriguing insight suggested by~(\ref{eqn:ce1}) is that
even when we choose to use hard thresholding with complete data,
if we adopt the self-consistency recipe~(\ref{eqn:selfcon}), we
should use ``soft" thresholding with incomplete data, as
(\ref{eqn:ce1}) is a form of soft thresholding. This can be seen
from the fact that as long as $\eta_l>0$, $0<|\tilde w_l^{(t)}|<
|w_l^{(t)}|$, as implied trivially by (\ref{eqn:estep}). It can
also been seen partially from the fact that as long as $\eta>0$,
$0<\beta(w, \eta)<1$, and therefore it provides an
``regression/shrinkage" effect. When $\eta \rightarrow 0$,
$\alpha(w, \eta) \rightarrow 0$ and $\beta(w, \eta) \rightarrow
1_{|w|\ge c}$, and thus (\ref{eqn:ce1}) goes back to the original
hard-thresholding, as it should be.

Fourth, although the updating expressions~(\ref{eqn:ce1})
and~(\ref{eqn:ce2}) for $\tilde w_{l}^{(t)}$ were derived for the
{\em hard} thresholding operation, it can be easily modified for
{\em soft} thresholding: $1_{(|w_l|\geq
c)}\mbox{sign}(w_l)\{|w_l|-c\}$. By similar calculations as in
Appendix~\ref{sec:closed}, for soft-thresholding, we only need to
add a simple term to $\tilde w_{l}^{(t)}$ of (\ref{eqn:ce1}),
relabeled as $\tilde w_{l, hard}^{(t)}$.  That is, we replace
(\ref{eqn:ce1}) by
\begin{equation}
\tilde w_{l, soft}^{(t)} = \tilde w_{l, hard}^{(t)} +
c\left\{\Phi\left(\frac{c-w_{l}^{(t)}}{\eta_l\sigma}\right)
-\Phi\left(\frac{c+w_{l}^{(t)}}{\eta_l\sigma}\right) \right\}.
\label{eqn:soft}
\end{equation}

Our numerical experiments suggest that the uses of the hard and
the soft thresholding operations inside the Ref algorithm give
similar practical performances, perhaps a reflection that both are
forms of ``soft" thresholding after all in the presence of
incomplete data. For this reason, for the rest of this paper we
shall concentrate on the hard thresholding operation in our
numerical experimentations.

\section{Modifications and Extensions}
\label{sec:extensions}
\subsection{Incorporating Smoothness Assumptions}
As we shall see from the simulation results in
Section~\ref{sec:sim}, when comparing to the interpolation based
methods, the self-consistency based procedures produce superior
estimates for $f(x_i)$ if $i\in\bI\obs$, but tend to produce
inferior estimates when $i\not\in\bI\obs$.  The reason for this is
as follows. Whenever an interpolation is employed to fill in the
missing responses $\{y_i:i\not\in\bI\obs\}$, some kind of
smoothness constraint is effectively imposed on
$\{f(x_i):i\not\in\bI\obs\}$. However, our basic self-consistency
procedures do not impose any prior smoothness assumptions on
$\{f(x_i):i\not\in\bI\obs\}$. Since many regression functions are
mostly smooth, it is reasonable to expect that methods that do not
impose any smoothness constraints on $\{f(x_i):i\not\in\bI\obs\}$
tend to produce inferior estimates for
$\{f(x_i):i\not\in\bI\obs\}$ than those methods that do impose
such constraints.  This is similar to comparing the maximum
likelihood procedure with a less efficient estimation procedure
but with a good constraint, or prior.  The latter can be better
than the former, not because it is better in retaining the
relevant information in the data, but because of the reasonable
constraint or prior information.

This analogy also suggests that, if such a smoothness assumption is
sensible, then this assumption should be included in the
self-consistent procedures; that is, to take advantage of both the
information available in the data and in the prior.  For example, if
linear interpolation is a good smoothing procedure to use, then the
MISC, Sim, and Ref algorithms can be further improved by modifying $\hat
f^{(t)}$ at the end of each iteration:
\begin{description}
\item [$\star$]
For each $i\not\in\bI\obs$, replace $\hat f^{(t)}(x_i)$ by the
following linearly interpolated value
\begin{equation}
\hat f^{(t)}(x_a) + \frac{\hat f^{(t)}(x_b) - \hat f^{(t)}(x_a)}
{x_b-x_a} (x_i-x_a),
\label{eqn:interp}
\end{equation}
where $x_a\in\bx\obs$ is the largest observed design point that is
less than $x_i$, and $x_b\in\bx\obs$ is the smallest observed design
point that is larger than $x_i$.
\end{description}

The above linear interpolation rule, of course, can be replaced by
many other interpolation rules if the corresponding smoothness
assumptions are more appropriate.  As demonstrated in
Section~\ref{sec:sim}, when such smoothness assumptions are
appropriate, these hybrid algorithms outperform both the pure
interpolation methods and the pure self-consistent methods in
terms of the mean squared errors integrated over both
$i\in\bI\obs$ and $i\not\in\bI\obs$.

\subsection{Non--Gaussian Errors}\label{sec:non}

As emphasized before, one of the main advantages of adopting the
self-consistency equation~(\ref{eqn:selfcon}) is that it is not
restricted to any particular model or error distribution, much
like the EM algorithm is not restricted to a particular class of
models. Indeed, the MISC algorithm is completely general, and can
be easily implemented to obtain estimators under non-Gaussian
errors: simply replace the Gaussian noise wavelet shrinkage
procedure in Step~2 of MISC with a suitable non--Gaussian
shrinkage method. A numerical demonstration with the Poisson model
will be given in Section~\ref{sec:poi}.

Given specific non-Gaussian thresholding rules, we can also derive
analytic approximations to the conditional expectation needed by
(\ref{eqn:selfcon}), under a specified non-Gaussian model. That
is, we can obtain refined algorithms for various other error
structures, including correlated errors. It is not possible to
give a general recipe here as how to perform such approximations,
as they need to be worked out on a case by case basis.  We emphasize
that this should not be viewed as a disadvantage of the
self-consistency method, because it is a general \emph{principle} for
constructing estimators, much like the implementation of the EM
algorithm always calls for individual derivations of its E-step (and
M-step) under the model of interest.

\subsection{Two--Dimensional Settings}\label{sec:2d}
In many imaging applications, especially in remote sensing area,
due to detector malfunction or some other reasons, the readings of
a small fraction of image pixels are missing.  These missing
values forbid the direct use of wavelet techniques for image
reconstruction.  This is an ideal setting for the self-consistent
procedures, as it is automatically an incomplete design problem
and the percentage of missing data tends to be small.

The generalizations of our methods from 1D settings to 2D settings
are in fact quite trivial.  Indeed, for all of the above
algorithms, the only modification needed is to replace the 1D DWT
with a 2D DWT. As for the 1D case, one would expect that the MISC
algorithm would give the best results, followed by Ref and then
Sim.  However, due to its computational cost, the 2D MISC
algorithm may not be practical.  Even for the 2D Ref algorithm,
the calculations of $\eta_l$ in (\ref{eqn:ce1})-(\ref{eqn:ce2})
can be lengthy, because one needs to calculate individual elements
of a 2D DWT matrix $\bW$. A very simple approximation is to
replace all $\eta_l$'s by their average, which is
\begin{equation}
\frac{1}{N}\mbox{trace}(\bI-\bW\bR\bW^\top) = 1 - \frac{n}{N}
\equiv C_m. \label{eqn:approxeta}
\end{equation}
This exceedingly simple approximation turned out working
surprisingly well in all our simulation studies, for both 1D and 2D
cases, and therefore we recommend its use whenever the more
refined calculations or approximations of $\eta_l$ are too
expensive.

Lastly we remark that in the image restoration or similar
contexts, the incorporation of an off-the-shelf interpolation step
often is not a good idea because real images tend to contain a
large amount of discontinuities (e.g., edges).

\section{Numerical Experiments}
\label{sec:sim} This section reports results of five sets of
numerical experiments that were conducted to study the empirical
properties of the above methods.  Throughout this section, the D5
wavelet of \citeasnoun{Daubechies92} was used as the mother
wavelet, and the primary resolution was 3.

\subsection{Visual Inspection}
\label{sec:vi1} Our first numerical experiment was simply to see
if the MISC (Section~\ref{sec:MISC}), the Sim
(Section~\ref{sec:Sim}) and the Ref (Section~\ref{sec:Ref})
algorithms, {\em without} using interpolation, would work at all.
We used the four well-known testing functions of
\citeasnoun{Donoho-Johnstone94}: {\em Heavisine}, {\em Doppler},
{\em Blocks} and {\em Bumps}.  For each of the test function, we
first simulated a regularly-spaced noisy data set as
in~(\ref{eqn:completed}), with $N=2048$ and signal--to--noise
ratio (snr) 7.  We then randomly deleted 30\% and 50\% of the
observations and applied the three algorithms to reconstruct the
curve, using the universal thresholding.

\begin{figure}
\begin{center}
\vspace*{-0.5cm}
\epsfig{file=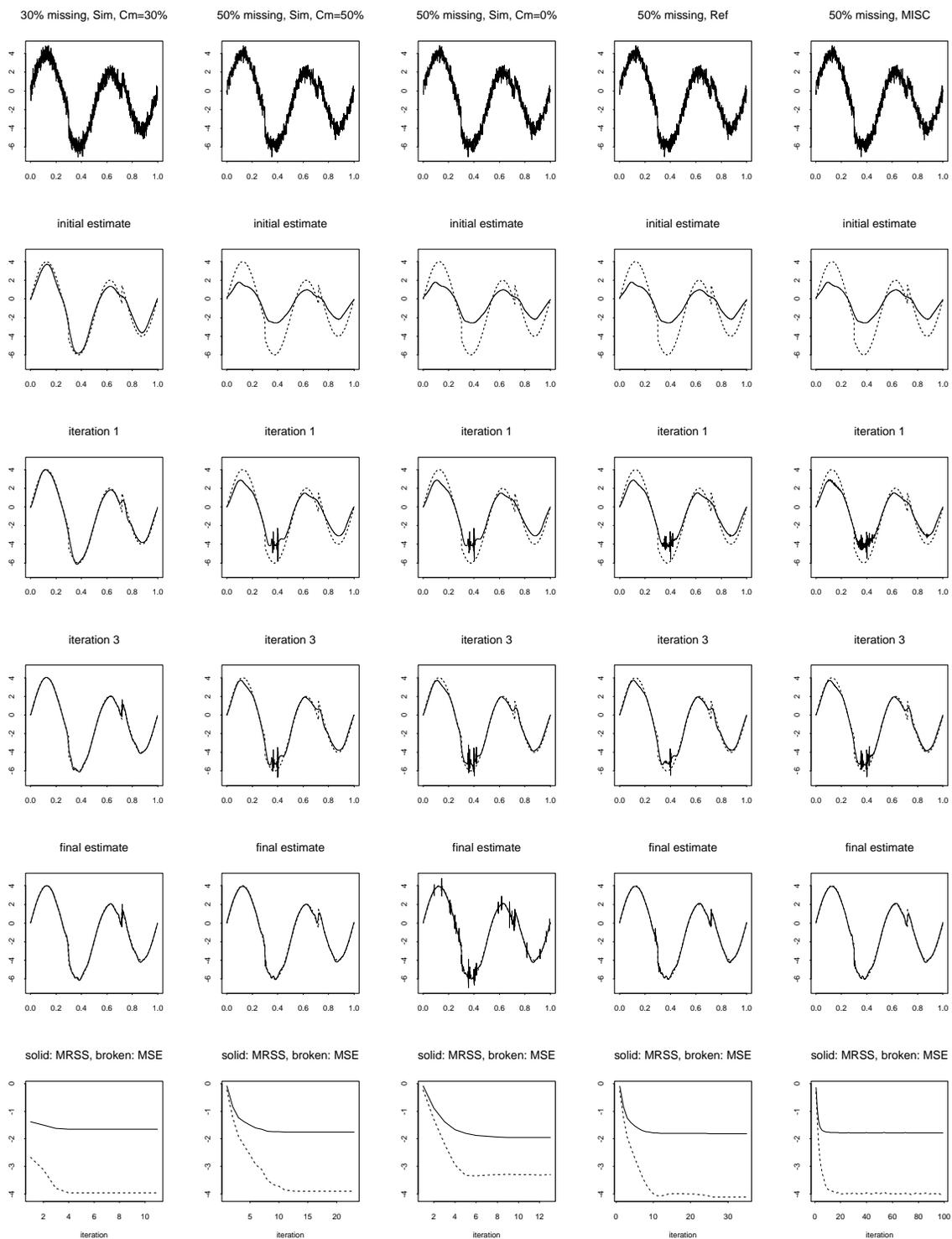,height=7.8in}
\vspace*{-0.5cm}
\end{center}
\caption{\small The performances of the Sim, Ref, and MISC algorithms for
recovering {\em Heavisine}.  For the second to the fourth rows, solid
curves are regression estimates and the dotted curve is the true
function.  See text for further details.}
\label{fig:ALLheavisine}
\end{figure}

The results for {\em Heavisine} are displayed in the first two
columns of Figure~\ref{fig:ALLheavisine} respectively. Each column
in Figure~\ref{fig:ALLheavisine} plots, from top to bottom, the
noisy incomplete data set; the initial estimate $\hat f^{(0)}$
obtained using the {\em S--Plus} function {\ttfamily lowess} with
a 10\% smoothing span; the first and the third iteration estimates
$\hat f^{(1)}$ and $\hat f^{(3)}$; the final estimate declared as
soon as reaching
$|\hat\sigma^{(t)}-\hat\sigma^{(t-1)}|/\hat\sigma^{(t)}<0.0001$;
and the plot of $\log(\mbox{MRSS}\obs)$ and $\log(\mbox{MSE}\obs)$
against the iteration step $t$.  Here $\mbox{MRSS}\obs$ is the mean
residual sum of squares and $\mbox{MSE}\obs$ is the mean squared
error calculated over all the observed design points: \bdm
\mbox{MRSS}\obs=\frac{1}{n}\sum_{i\in\bI\obs}\{y_i-\hat
f^{(t)}(x_i)\}^2 \quad \mbox{and} \quad
\mbox{MSE}\obs=\frac{1}{n}\sum_{i\in\bI\obs}\{f(x_i)-\hat
f^{(t)}(x_i)\}^2. \edm

Visually, we observe that the Sim algorithm produces good estimate
when $C_m=30\%$ (first column), and slightly worse estimate when
$C_m=50\%$ (second column). This worsening is expected because
with 50\% missing data there is a lot more uncertainty that could
be adequately captured by the variance inflation
formula~(\ref{eqn:inflate}). However, it still produces a much
better result than the naive approach that directly uses the
unadjusted $\tilde\sigma^2$ in Step~4, namely, by setting $C_m=0$
in (\ref{eqn:inflate}), as displayed in the third column. Results
of the MISC algorithm, with $M=100$, is displayed in the fifth
column for $C_m=50\%$.  This MISC algorithm does give a smaller
MSE than the Sim algorithm, but at an expense of increased
computational time.  It is because the MISC algorithm takes more
iterations and each
iteration is roughly $M=100$ times more expensive. A very exciting
observation is that the Ref algorithm provides essentially
identical results as MISC, as displayed in the fourth column, but
with a computational load (especially when approximating all
$\eta_l$'s by $C_m$) almost identical to that of the Sim
algorithm.

All the observations above were almost identically replicated with
the other three testing functions;  due to space limitation we
omit these displays.  For all these combinations, the Sim and the
Ref algorithms used $0.5$ to $3.5$ seconds user time on a Sun
Ultra 60 machine, depending on the testing function. For the MISC
algorithm, it took about $100\times T$ times longer, where $T$ is
the ratio of the number of iterations under MISC to that of Sim or
Ref ($T$ typically varies from 2 to 5).

\subsection{Effects of Interpolation and Approximations}
\label{sec:compint} Our second numerical experiment was conducted
to study the effects of the different approximations used in the
Sim and the Ref algorithms, using the MISC algorithm as a
benchmark for comparison.  The effects of interpolation is also
studied. The same four testing functions were used.  Other
experimental factors are: $N=512$, $\mbox{snr}=(5,7)$ and missing
percentage $C_m=(10\%, 30\%, 50\%)$.  For each combination of
testing function, snr and missing percentage, 100 incomplete data
sets were generated.  Then for each of these 100 generated data
sets, eight algorithms were applied to estimate $f$:
\begin{enumerate}
\item MISC, \item MISCI -- MISC with the interpolation
  step~(\ref{eqn:interp});
\item Sim, 
\item SimI --- Sim with~(\ref{eqn:interp});
\item Ref, 
\item RefI --- Ref with~(\ref{eqn:interp}); \item RefA --- Ref
with all $\eta_l$'s approximated by $C_m$, namely, by
(\ref{eqn:approxeta}), and \item RefAI
--- RefA with~(\ref{eqn:interp}).
\end{enumerate}
Instead of using the universal thresholding value
$\hat{\sigma}\sqrt{2\log N}$, in this experiment we follow
\citeasnoun{Antoniadis-Fan01} and used $\hat\sigma \sqrt{2\log N -
  \log(1+256\log N)}$.  \citeasnoun{Antoniadis-Fan01} showed that this
latter thresholding value is superior to the universal thresholding
value.

\begin{figure}[ht]
\begin{center}
\vspace*{-0.5cm}
\epsfig{file=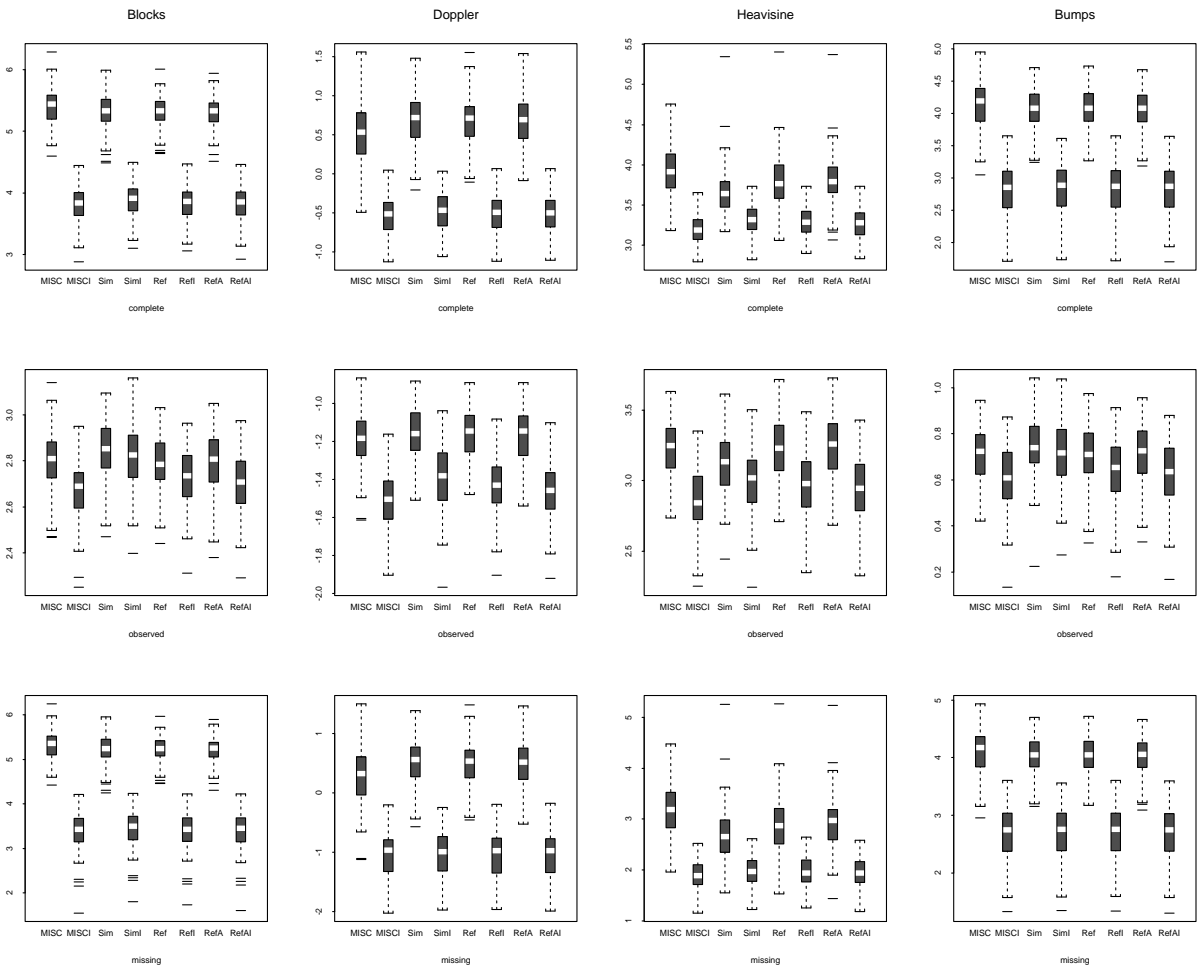,width=6.5in}
\vspace*{-1.5cm}
\end{center}
\caption{\small Boxplots of $\mbox{MSE}\com$, $\mbox{MSE}\obs$ and
  $\mbox{MSE}\mis$ for the eight algorithms tested in
  Section~\protect\ref{sec:compint}.}
\label{fig:boxplots1}
\end{figure}

Figure~\ref{fig:boxplots1} presents the boxplots, from top to bottom, of
$\mbox{MSE}\com$, $\mbox{MSE}\obs$, and $\mbox{MSE}\mis$, where each
column corresponds to a testing function.  Here $\mbox{MSE}\obs$ is the
same as in Section~\ref{sec:vi1}, and $\mbox{MSE}\mis$ and
$\mbox{MSE}\com$ are its counterparts summing over respectively all
the missing design points and all the design points.  The fraction of
missing data is $C_m=30\%$ and $\mbox{snr}=7$.  Results for
$C_m=(10\%,50\%)$ and $\mbox{snr}=5$ are similar and hence are omitted.

From the boxplots the following empirical observations can be made:
\begin{itemize}
\item Algorithms with interpolation are superior to their
  counterparts.  This illustrates the importance of employing suitable
  prior information in wavelet regression, and the above results
  helped to identify this importance by separating the efficiency
  inherited in the observed data from the prior information implicitly
  built into interpolation procedures.
\item When comparing results from Ref and RefA, and from RefI and RefAI,
  the $\eta$ approximation~(\ref{eqn:approxeta}) does not seem to have
  any adverse effects on Ref or RefI.
\item The performance of those Ref-type algorithms are very similar to
  the two MISC algorithms.  In fact for many experimental
  configurations, results from formal statistical tests (not reported
  here) suggest that the difference of the MSE values of RefAI and
  MISCI are statistically insignificant.
\end{itemize}
Given these observations, RefAI seems to be the best compromise, both
in terms of statistical performance and computational speed.

\subsection{Comparisons with Existing Methods}
\label{sec:compext} In this third experiment the algorithms SimI
and RefAI are compared to two popular wavelet regression
procedures for non--equispaced data. These two procedures are the
IRREGSURE procedure of \citeasnoun{Kovac-Silverman00} and the ROSE
procedure of \citeasnoun{Antoniadis-Fan01}.  The procedure
IRREGSURE is an interpolation based procedure where the SURE
method of \citeasnoun{Donoho-Johnstone95} is employed as the
thresholding procedure.  The ROSE procedure is a one-step
iterative procedure that uses penalized least squares and hard
thresholding.  Again, for SimI, RefAI and ROSE, the thresholding
value used was $\hat\sigma \sqrt{2\log N -
  \log(1+256\log N)}$.

\begin{figure}[ht]
\begin{center}
\vspace*{-0.5cm}
\epsfig{file=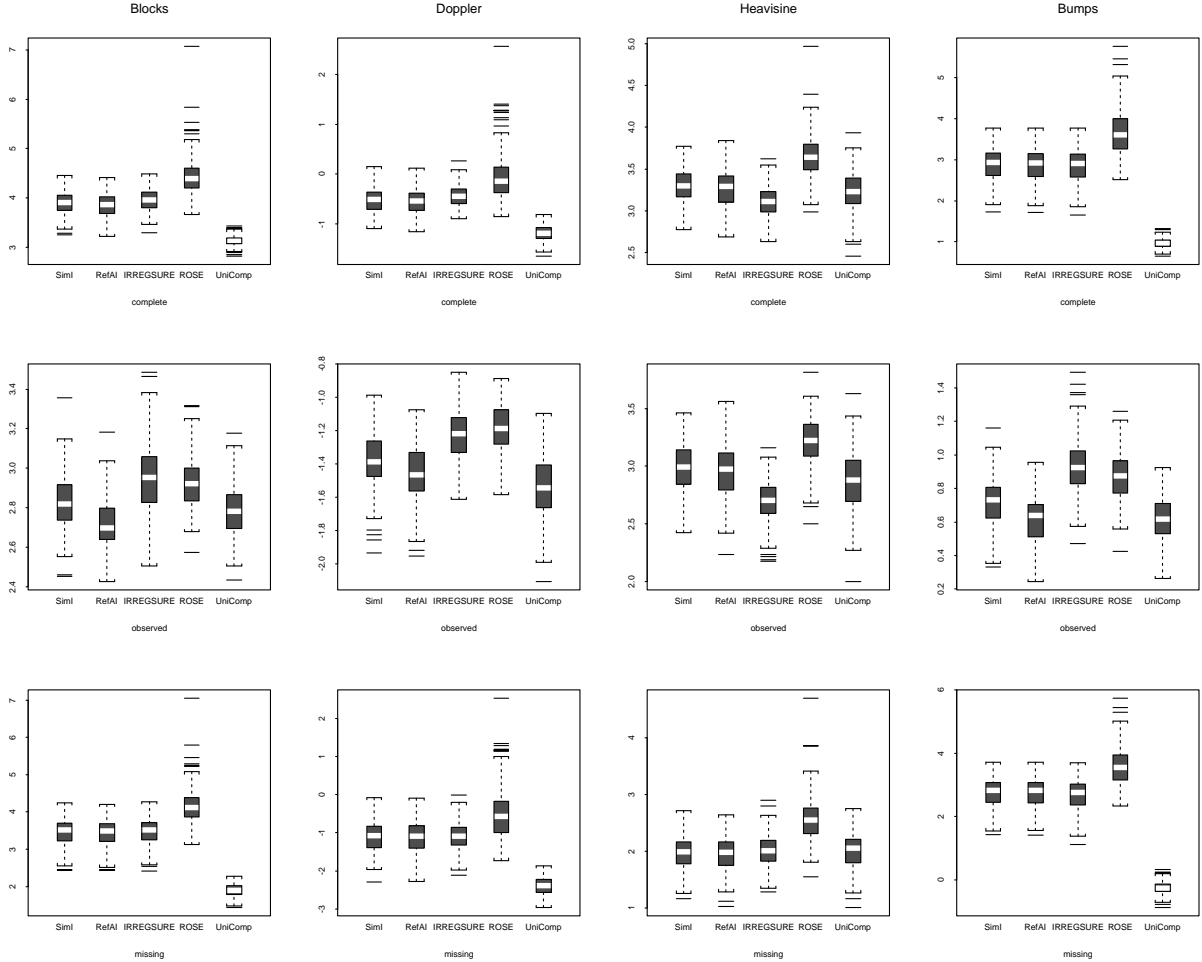,width=6.5in}
\vspace*{-1.5cm}
\end{center}
\caption{\small Boxplots of $\mbox{MSE}\com$, $\mbox{MSE}\obs$ and
  $\mbox{MSE}\mis$ for the four procedures compared in
  Section~\protect\ref{sec:compext}.}
\label{fig:boxplots2}
\end{figure}

We first present our results for those cases with $N=512$.  For each
combination of testing functions, $\mbox{snr}=(5,7)$ and $C_m=(10\%,
30\%, 50\%)$, 200 noisy data sets were generated.  For each noisy data
set, the four procedures SimI, RefAI, IRREGSURE and ROSE were applied to
obtain estimates for $f$, and their corresponding $\mbox{MSE}\com$,
$\mbox{MSE}\obs$, and $\mbox{MSE}\mis$ were computed.  In addition,
for the reason of providing a benchmark comparison, the universal
thresholding procedure with threshold $\hat\sigma \sqrt{2\log N -
  \log(1+256\log N)}$ was also applied to the {\em complete} noisy
data set.  We will label this procedure UniComp.  As the complete
data set was available to UniComp, it is expected that UniComp
would produce smaller MSE values than the other four procedures.
For the case $\mbox{snr}=7$ and $C_m=30\%$, boxplots of the
$\mbox{MSE}\com$, $\mbox{MSE}\obs$ and $\mbox{MSE}\mis$ values for
the five methods are displayed in Figure~\ref{fig:boxplots2} in a
similar fashion as before. Results under $\mbox{snr}=5$ and
$C_m=(10\%, 50\%)$ are similar and hence are omitted.

Pairwise Wilcoxon tests were applied to test if any two of the
four procedures have significantly different median values for
$\mbox{MSE}\com$, $\mbox{MSE}\obs$ and $\mbox{MSE}\mis$. The
significance level used was $1.25\%$ because of Bonferroni correction 
for multiple comparisons.
Based on the testing results, we ranked a procedure first if its
median MSE value is significantly less than those of the remaining
three, we ranked it second if its median MSE value is
significantly less than two but greater than the remaining one,
and so on so forth. If the median MSE values of two procedures are
not significantly different, they will share the same averaged
rank.  These rankings are listed in Table~\ref{table:ranking512}.
While no ranking method is perfect, such a ranking does provide a
good indicator of the relative merits of the methods being
compared.  Rankings under $\mbox{snr}=5$ are almost identical, and
thus are omitted.

From Figure~\ref{fig:boxplots2} and Table~\ref{table:ranking512}
one may conclude that RefAI is generally superior to the other
three procedures, SimI and IRREGSURE are roughly the same, while
ROSE is inferior.  A partial explanation for the relatively poorer
performance of ROSE is that it does not employ interpolation,
indicating potentially misleading comparative conclusions if we do
not distinguish between the information from the data and that
from (implicitly) imposed smoothness assumptions. By examining the
boxplots, one can see that, when comparing the $\mbox{MSE}\obs$
values, the performance of RefAI is in fact very similar to
UniComp. The above experiment was repeated for $N=2048$, but
without the ROSE procedure.  The relative rankings of RefAI, SimI
and IRREGSURE remain the same.

\renewcommand{\baselinestretch}{1.2}
\begin{table}[ht]
{\small
\begin{center}
\begin{tabular}{|cc|cccc|cccc|cccc|} \hline
& & \multicolumn{4}{c|}{MSE$\com$} & \multicolumn{4}{c|}{MSE$\obs$} &
 \multicolumn{4}{c|}{MSE$\mis$} \\ \cline{3-14}
function & {\bsy missing \% \esy}
 & \bsy{SimI}\esy & \bsy{RefAI}\esy & \bsy{IRREGSURE}\esy & \bsy{ROSE}\esy
 & \bsy{SimI}\esy & \bsy{RefAI}\esy & \bsy{IRREGSURE}\esy & \bsy{ROSE}\esy
 & \bsy{SimI}\esy & \bsy{RefAI}\esy & \bsy{IRREGSURE}\esy & \bsy{ROSE}\esy
\\ \hline
 & 10\% & 2 & 1 & 3 & 4 & 2 & 1 & 4 & 3 & 2 & 2 & 2 & 4 \\
{\em Blocks} & 30\% & 2 & 1 & 3 & 4 & 2 & 1 & 3.5 & 3.5 & 2.5 & 1 & 2.5 & 4 \\
 & 50\% & 2 & 1 & 3 & 4 & 2 & 1 & 3 & 4 & 2.5 & 1 & 2.5 & 4 \\ \hline
 & 10\% & 2 & 1 & 3.5 & 3.5 & 2 & 1 & 4 & 3 & 1.5 & 1.5 & 3 & 4 \\
{\em Doppler} & 30\% & 2 & 1 & 3 & 4 & 2 & 1 & 3 & 4 & 2 & 2 & 2 & 4 \\
 & 50\% & 3 & 1.5 & 1.5 & 4 & 2 & 1 & 3 & 4 & 3 & 2 & 1 & 4 \\ \hline
 & 10\% & 2 & 1 & 4 & 3 & 2 & 1 & 4 & 3 & 1.5 & 1.5 & 4 & 3 \\
{\em Heavisine} & 30\% & 3 & 2 & 1 & 4 & 3 & 2 & 1 & 4 & 2 & 2 & 2 & 4 \\
 & 50\% & 2.5 & 2.5 & 1 & 4 & 2 & 3 & 1 & 4 & 2.5 & 2.5 & 1 & 4 \\ \hline
 & 10\% & 2 & 1 & 3 & 4 & 2 & 1 & 4 & 3 & 2.5 & 2.5 & 1 & 4 \\
{\em Bumps} & 30\% & 3 & 1.5 & 1.5 & 4 & 2 & 1 & 4 & 3 & 3 & 2 & 1 & 4 \\
 & 50\% & 3 & 2 & 1 & 4 & 2 & 1 & 3 & 4 & 3 & 2 & 1 & 4 \\ \hline
\multicolumn{2}{|c|}{average rank} & 2.4 & 1.4 & 2.4 & 3.8 & 2.1 & 1.3
 & 3.1 & 3.5 & 2.4 & 1.8 & 1.9 & 3.9 \\ \hline
\end{tabular}
\end{center}
}
\caption{\small Wilcoxon rankings for the four wavelet regression
  procedures compared in Section~\ref{sec:compint} when $N=512$.}
\label{table:ranking512}
\end{table}
\renewcommand{\baselinestretch}{1.69}

\subsection{Poisson Model}\label{sec:poi}

\begin{figure}[ht]
\vspace*{-1cm}
\begin{center}
\epsfig{file=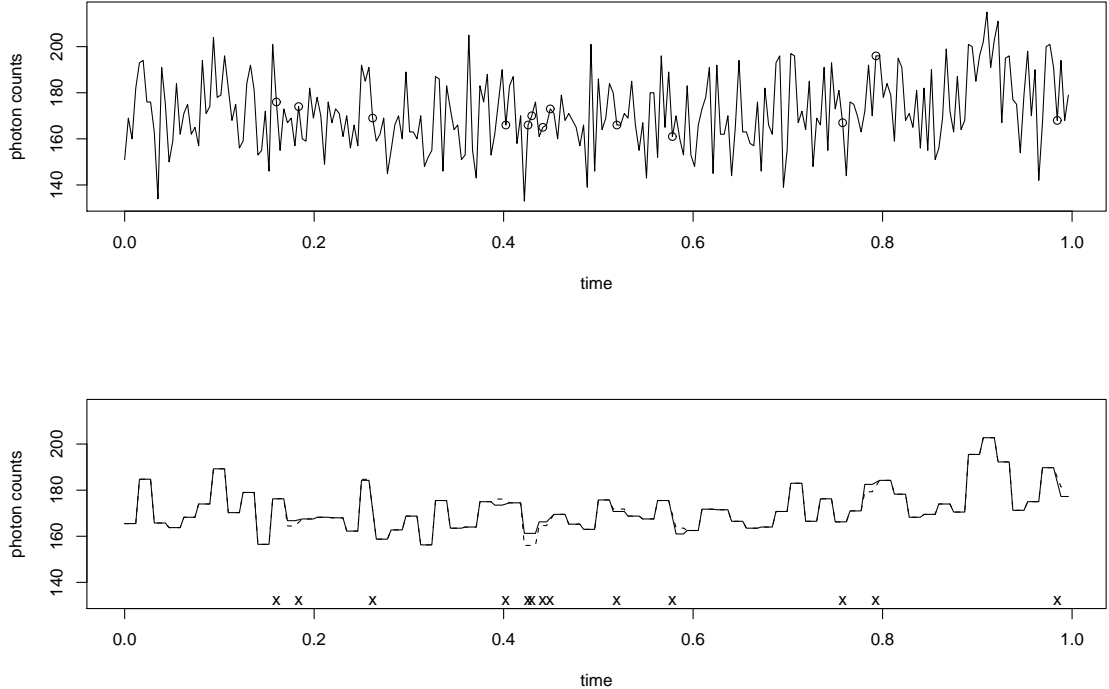,height=6in,angle=270}
\vspace*{-0.5cm}
\end{center}
\caption{\small Top panel: photon counts from the collapsed star
RXJ1856.5-3754.  Circles indicate (artificial) missing
observations.
  Bottom panel: reconstructed curves using the complete data (solid
  line) and missing data (broken line).  The crosses mark the
  locations of the missing data.}
\label{fig:poisson}
\end{figure}

Displayed in the top panel of Figure~\ref{fig:poisson} are the
Poisson photon counts, captured at 256 time intervals, from the
collapsed star RXJ1856.5-3754,  which is about 400 light years
from Earth in the constellation Corona Australis.  Note that in
this plot the time index has been re-scaled to $[0,1]$.  One
reason that astrophysics scientists are interested in this
collapsed star is that they believe its matter is even denser than
nuclear matter, the most dense matter found on Earth. We refer
interested readers to {\tt
http://chandra.harvard.edu/photo/2002/0211/index.html} for the
scientific issues surrounding this data set.

From this Poisson data set the following two smoothed curve
estimates are obtained.  The first one was constructed from the
complete data set while the second one was constructed with the
presence of missing data.  The complete data reconstructed curve
was obtained by applying the Poisson wavelet thresholding rule of
\citeasnoun[Equation~(13)]{Kolaczyk99b} to the full data set, and
it is displayed as the solid line in the bottom panel of
Figure~\ref{fig:poisson}. The missing-data curve estimate was
obtained as follows.  Firstly 5\% of the data points were removed.
These missing data points are marked by ``{\tt o}'' in the top
panel and their locations are indicated by ``{\tt x}'' in the
bottom panel of Figure~\ref{fig:poisson}.  Then the MISC algorithm
was applied to this missing data set with the same Poisson
thresholding rule and $M=100$.  The resulting curve estimate is
the broken line in the bottom panel of Figure~\ref{fig:poisson}.

From these plots one could notice that in regions with no missing
data (e.g., for $t$ in $[0, 0.15]$ and $[0.6,0.75]$) or when the
values of the missing data are not local extrema (e.g., at
$t\approx 0.26$ and $t\approx 0.76$), the complete-data and the
missing-data estimates are virtually the same.  However, the two
estimates do have small differences at regions where missing data
are clustered (e.g., at $t\approx 0.44$) or when the values of the
missing data are local extrema (e.g., $t\approx 0.19$ and
$t\approx 0.79$). This simple example therefore illustrates both
the feasibility of the MISC algorithm for non-Gaussian data, and
the fact that wavelets methods have the ability to localize the
damage caused by the incomplete observations.

\subsection{Image Denoising}
\label{sec:image} In this last experiment we explore the
performance of our algorithms in the context of image denoising.
We are only aware of a very limited number of existing methods
that are specifically designed to perform image denoising with
missing data. One method was described by
\citeasnoun{Naveau-Oh04}, which although can be applied to handle
missing values around the image edges, was primarily proposed to
reduce boundary effects. Their updating algorithm is similar to
our Sim algorithm, but without the variance inflation
formula~(\ref{eqn:inflate}), and therefore inferior results are
expected. The other is by \citeasnoun{Hirakawa-Meng06} using an
EM-type approach for simultaneous demosicing and denosing, which in
fact was motivated by our current work.

\begin{figure}[ht]
\begin{center}
\begin{tabular}{cc}
\epsfig{figure=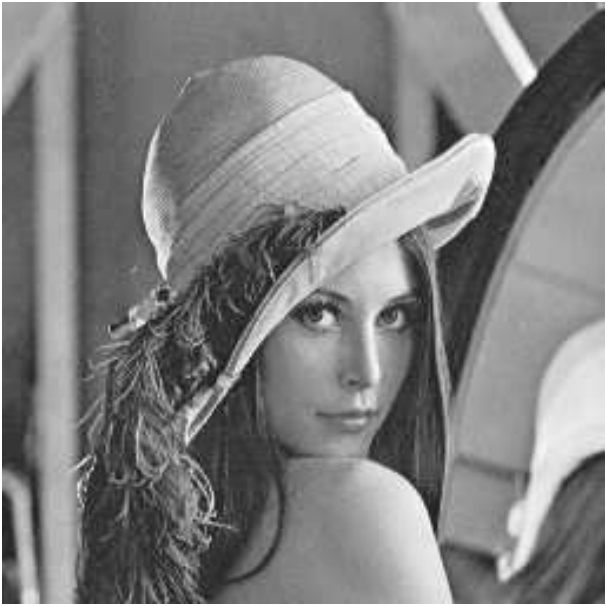,width=8.0cm} &
\epsfig{figure=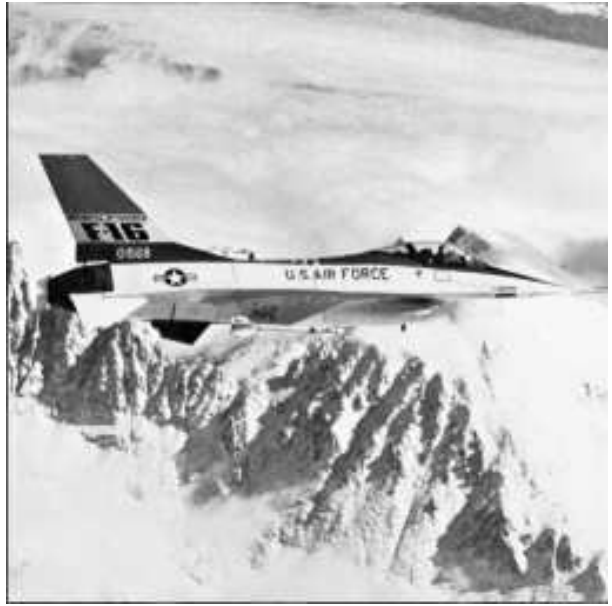,width=8.0cm} \\
(a) Lena & (b) Airplane
\end{tabular}
\end{center}
\caption{\small Testing images used in the image denoising experiment.}
\label{fig:testimages}
\end{figure}

Three 2D algorithms were studied: the MISC (with $M=10$), the Sim,
the RefA procedures. Two testing images of size $256\times 256$
were used: the well-known Lena image displayed in
Figure~\ref{fig:testimages}(a) and the Airplane image displayed in
Figure~\ref{fig:testimages}(b). Also, two snrs and three missing
data percentages were tested: $\mbox{snr}=(5, 7)$ and $C_m=(10\%,
30\%, 50\%)$.  Lastly, two missing data formation mechanisms were
tested.  The first one is missing at random, in which missing
pixel locations were randomly selected from the image, while in
the second mechanisms the missing pixels were clustered together.
Note that because of the computational cost, it was too costly to
run MISC with $M=100$ for our simulation studies, which typically
takes roughly 1 hour for one replicate on a Sun Ultra 60 machine.

For each of the above experimental factor combinations, 100 noisy
images were generated, and the above three algorithms were applied
to reconstruct the corresponding true images, using the adjusted
universal thresholding value: $\hat\sigma \sqrt{2\log N -
\log(1+256\log N)}$. As to provide a benchmark for comparison, for
each noisy image, we also applied the universal denoising method
\cite{Donoho-Johnstone94}, with the same adjusted thresholding
value,  to reconstruct the corresponding true image using the
complete data.  As before, we refer this method as UniComp.  As
UniComp has the full information from $\by$, it is expected that
it would produce better reconstructed images than the other three
algorithms.

\begin{figure}[ht]
\centerline{\epsfig{figure=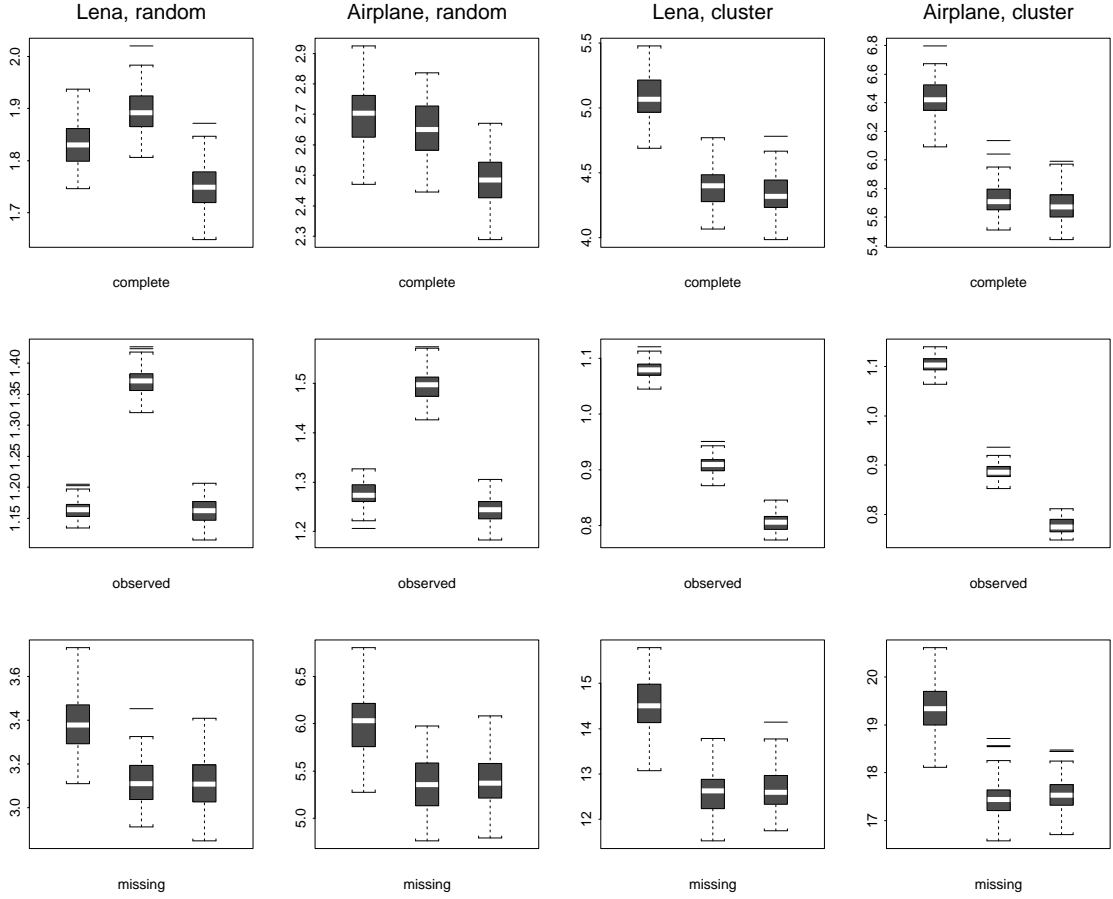,width=6in}}
\caption{\small Boxplots of the MSE ratios resulted from the image
  denoising
  experiment in Section~\protect\ref{sec:sim}.  In each panel the
  left, middle and right boxplots correspond, respectively, to the
  MISC, Sim and RefA algorithms.}
\label{fig:imageboxplots}
\end{figure}

For every reconstructed images, we calculated $\mbox{MSE}\com$,
$\mbox{MSE}\obs$ and $\mbox{MSE}\mis$ as measures of
reconstruction quality. (We are aware of the fact that MSE is not
a good measure for visual quality, but in the absence of a
commonly agreed measure for visual quality, the MSE still serves
as a statistically useful criterion for comparisons). In addition
we also computed the following MSE ratio: \bdm
r\com(\mbox{MISC})=\frac{\mbox{$\mbox{MSE}\com$ of MISC}}
                      {\mbox{$\mbox{MSE}\com$ of UniComp}}.
\edm
Similar MSE ratios for the observed ($r\obs(\mbox{MISC})$) and missing
data ($r\mis(\mbox{MISC})$), and for the Sim and RefA algorithms
($r\com(\mbox{Sim})$, $r\obs(\mbox{Sim})$,
$r\mis(\mbox{Sim})$, $r\com(\mbox{RefA})$,
$r\obs(\mbox{RefA})$ and $r\mis(\mbox{RefA})$)
were also calculated.  Since UniComp reconstructed the images with the
complete data, it is expected that all these MSE ratios are bigger
than 1.  For $\mbox{snr}=7$ and $C_m=30\%$, boxplots of these MSE
ratios are given in Figure~\ref{fig:imageboxplots}.  Boxplots for
other experimental settings are similar and hence omitted.  From
Figure~\ref{fig:imageboxplots} some major empirical conclusions can be
obtained.

First, as all $r\obs(\mbox{MISC})$ values are fairly close to 1,
the easy-to-implement benchmark MISC algorithm performs reasonably
well for those observed pixels.  Secondly, the RefA algorithm is
superior to the other two algorithms, as it does not require
multiple imputation (as opposed to MISC) and it uses a better
approximation than the Sim algorithm.  Lastly, an unexpected
observation is that, $r\obs(\mbox{RefA})$ is in fact less than 1
when the locations of the missing data are clustered together. We
currently do not have an explanation for this phenomenon, other
then noting that hidden biases resulting from model defects can be
more pronounced with more data.

\begin{figure}[ht]
\begin{center}
\begin{tabular}{cc}
\epsfig{figure=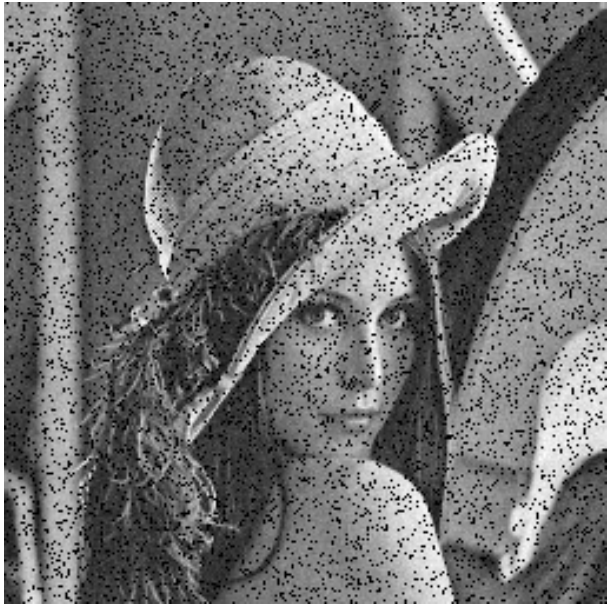,width=8.0cm} &
\epsfig{figure=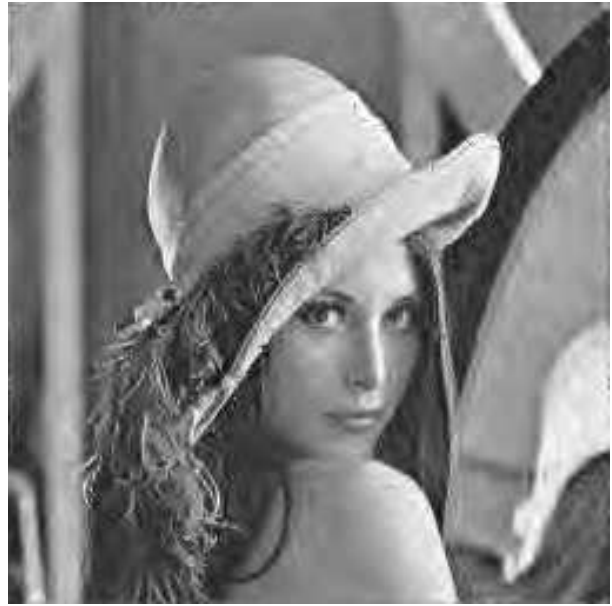,width=8.0cm} \\
(a) Degraded Lena & (b) Reconstructed Lena
\end{tabular}
\end{center}
\caption{\small Degraded (a) and reconstructed (b) Lena when the
  pixels are missing at random.}
\label{fig:lenarandom}
\end{figure}

\begin{figure}[ht]
\begin{center}
\begin{tabular}{cc}
\epsfig{figure=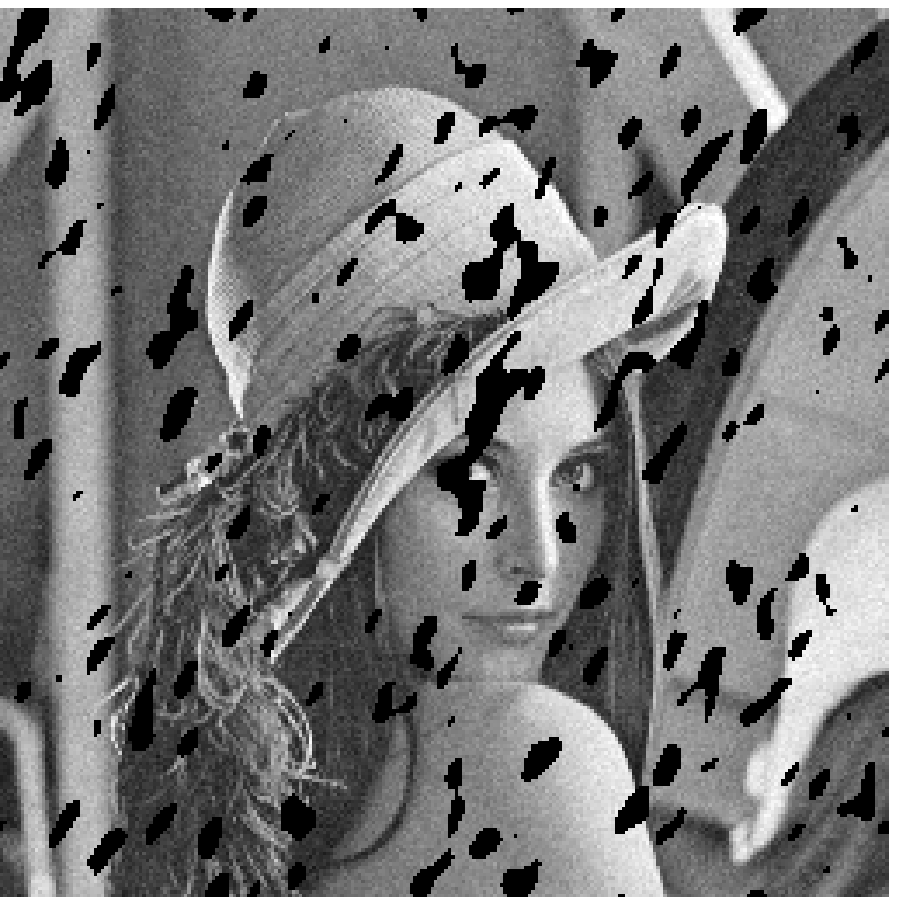,width=8.0cm} &
\epsfig{figure=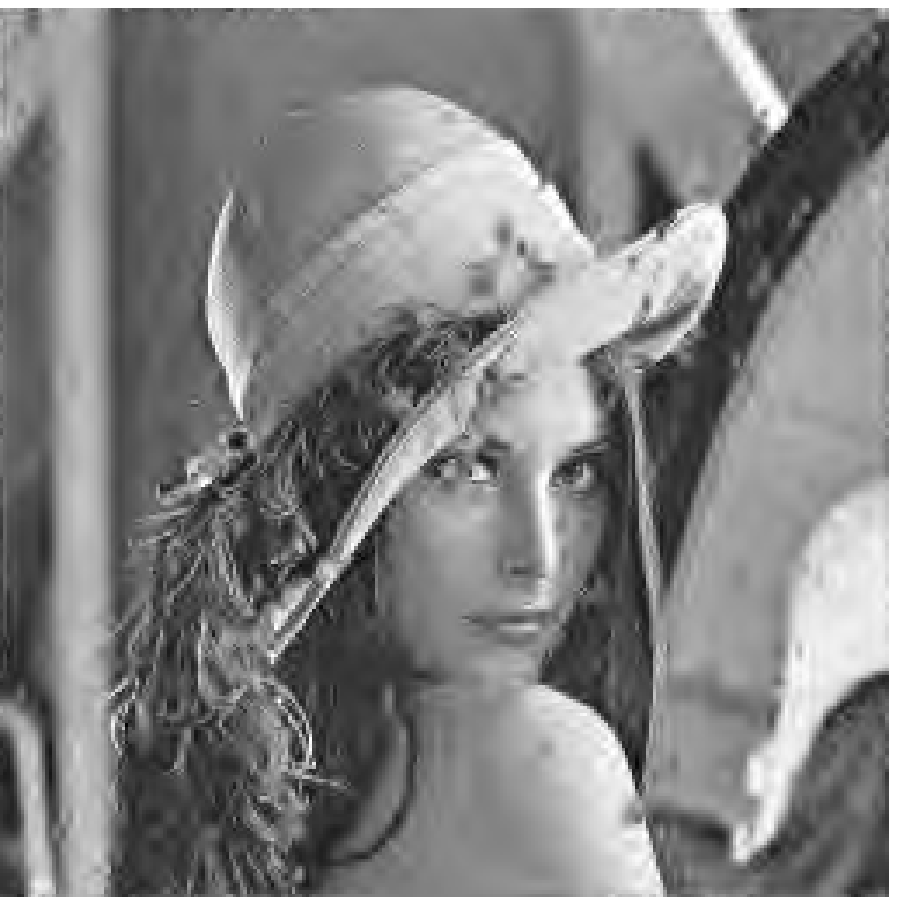,width=8.0cm} \\
(a) Degraded Lena & (b) Reconstructed Lena
\end{tabular}
\end{center}
\caption{\small Similar to Figure~\protect\ref{fig:lenarandom} but for
  clustered missing pixels.}
\label{fig:lenacluster}
\end{figure}

For the purpose of visual inspection, two degraded versions of
Lena are displayed in Figures~\ref{fig:lenarandom}(a)
and~\ref{fig:lenacluster}(a).  Those black pixels represent the
locations of the missing values.  The snr is 7 and the missing
percentage is 10\%.  Figures~\ref{fig:lenarandom}(b)
and~\ref{fig:lenacluster}(b) display the corresponding
reconstructed images obtained from the RefA algorithm. The quality
of the reconstructed ones is quite acceptable. The one with
clustered missing data is particularly impressive, especially
considering that the method we used did not take into account the
cluster nature of the missing data. Reconstruction algorithms as
such are particularly useful for image inpainting (e.g., see
\myncite{Criminisi-et-al04} and \myncite{Tschumperle-Deriche05}).

\section{Summary and Future Works}
\label{sec:conclude} A main goal of this paper is to demonstrate that
the self-consistency principle is a very versatile and fruitful method
for dealing with wavelet modeling, and more generally with
non-parametric and semi-parametric regressions, when facing incomplete
data.  By viewing irregularly-spaced data as a form of incomplete
data, it also provides a rather general methodology for wavelet
reconstructions with irregularly-spaced data by taking advantage of
the existence of those well-studied methods developed for
regularly-spaced data, much in the same spirit as with the EM 
algorithm or multiple imputation. Indeed, the specific algorithms
we proposed here directly use either multiple imputation or steps
very similar to the E-step (and M-step) of the EM algorithm.

Much remains to be done, of course. The most urgent ones are
theoretical properties of the estimators and algorithms we
propose. Simulations were crucial in our development of the
estimators and algorithms given in this paper, but they are no
substitute of rigorous theoretical investigations. It is therefore
our hope that our empirical findings are convincing, or at least
suggestive, enough that they would motivate a general theoretical
investigation of the self-consistent wavelet estimators, or more
generally self-consistent regression estimators. It should also be
of great interest to investigate the theoretical connections
between the self-consistent wavelet estimators and other
constructions of wavelet estimators with irregularly-spaced
data, such as via lifting (e.g., \myncite{Delouille-et-al04} and
\myncite{Nunes-et-al06}).

On the methodological side, besides developing even more refined
algorithms, especially for non-Gaussian and correlated errors, a
logical next step is to combined the self-consistent principle
with Bayesian methods. Indeed, with Bayesian methods, the
dealing with missing data is done jointly with the inference of
parameters and regression functions. Preliminary work has shown
great promise, as reported in \citeasnoun{Hirakawa-Meng06}.  We
are currently investigate the self-consistent approach with a
number of Bayesian methods for wavelet reconstructions, including
over-complete expansions (e.g.,
\myncite{Abramovich-et-al98,Johnstone-Silverman05} and
\myncite{elad06:_stabl}).

\section*{Acknowledgment}
The authors thank the National Science Foundation for partial
support (Grants No.~0203901 and No.~0204552) for this work,
particularly for the funding to present an abbreviated and
preliminary version (\myncite{Lee-Meng05}) at the 30th IEEE
International Conference on Acoustics, Speech, and Signal
Processing.

\appendix
\section{Derivations of (\ref{eqn:ce1}), (\ref{eqn:ce2}) and
  (\ref{eqn:soft})}
\label{sec:closed} To show (\ref{eqn:ce1}) and (\ref{eqn:ce2}),
let $w_l \sim \N(d, \tau^2)$ and write $Z=\frac{w_l-d}{\tau}$;
i.e., $Z$ is standard normal. We will also need the following
fact. Let $Z \sim \N(0,1)$, and denote its probability density
function by $\phi(z)$. Then, for any constant $c$,
$E\{1_{(Z>c)}Z\}=\phi(c)$ and $E\{1_{(Z<-c)}Z\}=-\phi(c)$.

With this setup, we have
\begin{eqnarray*}
E\{1_{(|w_l|\geq c)}w_l\}
& = &
E\{1_{(d+\tau Z\geq c)}(d+\tau Z)\} + E\{1_{(d+\tau Z\leq -c)}(d+\tau Z)\} \\
& = &
d P\left(Z\geq \frac{c-d}{\tau}\right) + \tau E\{1_{(Z\geq
  \frac{c-d}{\tau})}Z\} +
d P\left(Z\leq \frac{-(c+d)}{\tau}\right) + \tau E\{1_{(Z\leq
  \frac{-(c+d)}{\tau})}Z\} \\
& = &
d \left\{1-\Phi\left(\frac{c-d}{\tau}\right)\right\} +
\tau\phi\left(\frac{c-d}{\tau}\right) +
d \left\{1-\Phi\left(\frac{c+d}{\tau}\right)\right\} -
\tau\phi\left(\frac{c+d}{\tau}\right) \\
& = &
d \left\{2-\Phi\left(\frac{c-d}{\tau}\right) -
           \Phi\left(\frac{c+d}{\tau}\right)\right\} +
\tau\left\{\phi\left(\frac{c-d}{\tau}\right) -
  \phi\left(\frac{c+d}{\tau}\right) \right\}.
\end{eqnarray*}
Now by substituting $d=w_l^{(t)}$ and $\tau=\eta_l\sigma$ into the
above expression, we obtain~(\ref{eqn:ce1}) and~(\ref{eqn:ce2}).

Equation~(\ref{eqn:soft}) is derived using essentially the same
steps as above, except the whole calculation begins with
$E[1_{(|w_l|\geq c)}\mbox{sign}(w_l)\{|w_l|-c\}]$, which only
differs from the hard-thresholding formula by the simple term $-c
E\{1_{(|w_l|\geq c)}\mbox{sign}(w_l)\}= -c \{P(w_l\ge c) -
P(w_l\le -c)\}$.  This can be easily seen to be corresponding to the 
second term on the right-hand side of~(\ref{eqn:soft}).

\bibliography{refall,meng}
\bibliographystyle{myagsm}
\end{document}